\documentclass[a4paper]{amsart}
\usepackage{pdfsync}

\copyrightinfo{2006}
  {O.Lakkis \& Ch.Makridakis}

\PII{}

         \usepackage{ifthen}
         \usepackage{amsmath}
         \usepackage{amssymb}    
         \usepackage{stmaryrd} 
         \usepackage{eufrak}     
         \usepackage{mathrsfs} 
         \usepackage{bbold}    
         \provideboolean{usemathrsfs}
         \setboolean{usemathrsfs}{true}
         \usepackage{enumerate}
         \usepackage{xspace}
         \usepackage[greek,british,english]{babel}
         \usepackage{boxedminipage}
         \usepackage[utf8]{inputenc}
         \usepackage{verbatim}
         \usepackage{fancyvrb}
         \usepackage{listings}
         \usepackage{hyperref}
         \provideboolean{usebeamer}
         \provideboolean{isthesis}
         \provideboolean{isamsltex}
         \provideboolean{issiamltex}
     \usepackage{hyperref}
     
     \providecommand{\linkedemail}[1]{{\makeatletter\color{blue}\href{mailto:#1}{#1}\makeatother}}
     \providecommand{\Ignore}[1]{}
     \providecommand{\freeze}[1]{}
     \providecommand{\crossout}[1]{{\textcolor{red!20}{#1}}}
     \providecommand{\highlight}[1]{{\color{blue}#1}}
     \providecommand{\standout}[1]{\colorbox{a}{\textcolor{g}{#1}}}
     \provideboolean{shownotes}
     \setboolean{shownotes}{true}
     \newcounter{margnote}[page]
     \providecommand{\margnotemark}{{\standout{\upshape\texttt{\arabic{margnote}}}}}
     \providecommand{\margnote}[2][]{
       \ifthenelse{
         \boolean{shownotes}
       }{
         \stepcounter{margnote}
         \margnotemark
         \marginpar{
           \texttt{
             \begin{minipage}{2cm}
               \raggedright\tiny
               \margnotemark{#1}: 
               #2
             \end{minipage}
           }
         }
       }{
       }
     }
     \providecommand{\mathnote}[2][]{
       \ifthenelse{
         \boolean{shownotes}
       }{
         \stepcounter{margnote}
         \margnotemark
         \text{
           \standout{
             \texttt{
               \tiny
               #1:
               #2
             }
           }
         }
       }{
       }
     }
     \provideboolean{showtodo}

     \providecommand{\todo}[1]{\ifthenelse{\boolean{showtodo}}{\margnote[To do.]{#1}}{}}
     \providecommand{\Todo}[1]{
       \ifthenelse{\boolean{showtodo}}{
         \setlength{\parindent}{0pt}\ 
         \par
         \begin{boxedminipage}{\textwidth}
           \texttt{#1}
         \end{boxedminipage}
         \par
       }{}}

     \provideboolean{showcomments}
     \providecommand{\margincomment}[1]{
     \ifthenelse{\boolean{showcomments}}{\marginpar{\tiny #1}}{}
     }
     \provideboolean{showchanges}
     \setboolean{showchanges}{false}
     \providecommand{\changes}[1]{
       \ifthenelse{\boolean{showchanges}}{{\highlight{#1}}}{#1}
     }
     \providecommand{\changefromto}[3][replace with]{
       \ifthenelse{\boolean{showchanges}}
       {{\crossout{#2}\margnote{#1}}{\highlight{#3}}}
       {#3\xspace}
     }
     \providecommand{\ChangePar}[2]{
       \ifthenelse{\boolean{showchanges}}
       {{\par$\mapsfrom$ \textcolor{red!20}{#1}}{\par$\mapsto$ \textcolor{blue}{#2}}}
       {\par #2}
     }
     \providecommand{\InsertPar}[1]{
       \ifthenelse{\boolean{showchanges}}
       {{\par$\mapsto$ \textcolor{blue}{#1}}}
       {\par #1}
     }

     \providecommand{\ie}{\ensuremath{\text{ i.e., }}\xspace}

     \providecommand{\aposteriori}{{aposteriori}\xspace}
     \providecommand{\Aposteriori}{{Aposteriori}\xspace}
     
     \usepackage{mathrsfs}
     \provideboolean{usemathrsfs}
     \setboolean{usemathrsfs}{true}
     \providecommand{\mathscript}
     	   {\mathscr}
     
      \providecommand{\cA}{\ensuremath{\mathscript A}\xspace}

      \providecommand{\cE}{\ensuremath{\mathscript E}\xspace}

      \providecommand{\cI}{\ensuremath{\mathscript I}\xspace}

      \providecommand{\cM}{\ensuremath{\mathscript M}\xspace}

      \providecommand{\cR}{\ensuremath{\mathscript R}\xspace}
      \providecommand{\cS}{\ensuremath{\mathscript S}\xspace}
      \providecommand{\cT}{\ensuremath{\mathscript T}\xspace}
      
      \providecommand{\cV}{\ensuremath{\mathscript V}\xspace}

      \providecommand{\rN}{\ensuremath{\mathbb N}\xspace}
      
      \providecommand{\rP}{\ensuremath{\mathbb P}\xspace}
      
      \providecommand{\rR}{\ensuremath{\mathbb R}\xspace}

      \providecommand{\rV}{\ensuremath{\mathbb V}\xspace}
      \providecommand{\rW}{\ensuremath{\mathbb W}\xspace}

      \providecommand{\rZ}{\ensuremath{\mathbb Z}\xspace}
     
      \providecommand{\naturals}{\rN\xspace}
      \providecommand{\Nat}{\ensuremath{\rZ^+}}

      \providecommand{\reals}{\rR}
      \providecommand{\R}[1]{\reals^{#1}}

      \providecommand{\one}{\ensuremath{\mathbb 1}\xspace}

      \providecommand{\Tolto}{\smallsetminus}

      \providecommand{\inner}{\cdot}

      \providecommand{\W}{\ensuremath{\varOmega}\xspace}

      \providecommand{\ep}{\ensuremath{\varepsilon}\xspace}
      \providecommand{\epsi}{\ensuremath{\epsilon}\xspace}

      \providecommand{\w}{\ensuremath{\omega}\xspace}
     
      \providecommand{\qp}[1]{\ensuremath{\!\left({#1}\right)}}
      \providecommand{\qpreg}[1]{\ensuremath{(#1)}}

      \providecommand{\qpbigg}[1]{\ensuremath{\bigg(#1\bigg)}}
      \providecommand{\qpBigg}[1]{\ensuremath{\Bigg(#1\Bigg)}}
      \providecommand{\qb}[1]{\ensuremath{\!\left[{#1}\right]}}

      \providecommand{\qc}[1]{\ensuremath{\left\{{#1}\right\}}}

      \providecommand{\qa}[1]{\ensuremath{\left\langle{#1}\right\rangle}}

      \providecommand{\opinter}[2]{\ensuremath{\left(#1,#2\right)}\xspace}
      \providecommand{\clinter}[2]{\ensuremath{\left[#1,#2\right]}\xspace}
      \providecommand{\opclinter}[2]{\ensuremath{\left(#1,#2\right]}\xspace}
      
      \providecommand{\powqp}[2]{\ensuremath{\qp{#2}^{\kern -.1em\lower .2ex\hbox{\scriptsize $#1$}}\kern-.1em}}
      \providecommand{\powqpreg}[2]{\ensuremath{\qpreg{#2}^{\kern -.2em\lower .3ex\hbox{\scriptsize $#1$}}\kern-.3em}}
      \providecommand{\psqrt}[1]{\powqp{1/2}{#1}}
      
      \providecommand{\norm}[1]{\ensuremath{\left|#1\right|}}

      \providecommand{\Norm}[1]{\ensuremath{\left\|#1\right\|}}

      \providecommand{\ltwop}[2]{\ensuremath{\qa{#1,#2}}}

      \providecommand{\duality}[2]{\ensuremath{\left\langle #1\,\vert\,#2\right\rangle}}
      
      \providecommand{\ensemble}[2]{\ensuremath{\left\{ #1:\;#2 \right\}}}

      \providecommand{\setof}[1]{\qc{#1}}

      \providecommand{\seqof}[1]{\qp{#1}}
      \providecommand{\seq}[1]{\seqof{#1}}
      \providecommand{\seqi}[2]{\seq{#1_{#2}}}
      \providecommand{\seqs}[2]{\seq{#1}_{#2}}

      \providecommand{\sumifromto}[3]{\ensuremath{\sum_{#1=#2}^{#3}}}

      \providecommand{\jump}[1]{\ensuremath{\left\llbracket #1\right\rrbracket}}
     
      \providecommand{\fromto}[2]{\ensuremath{\left[#1:#2\right]}}

      \providecommand{\integerbetween}[2]{\ensuremath{\in\fromto{#1}{#2}}}
      \providecommand{\rangefromto}[3]{\ensuremath{#1\integerbetween{#2}{#3}}}
     
      \renewcommand{\d}{\ensuremath{\,\mathrm{d}}}
      \providecommand{\ds}[1]{\ensuremath{\,\mathrm{s}(\mathrm d#1)}}
      \providecommand{\D}{\ensuremath{\,\mathrm{D}}}
     
     \providecommand{\registered}%
         {\ensuremath{{}^{\bigcirc\!\;\!\!\!\!\!\!\!\;\text{\sc r}}}}
     
     \providecommand{\AND}{\ensuremath{\text{ and }}}

     \providecommand{\ball}{\operatorname{B}}                  
     \providecommand{\Ball}[2]{\ball_{{#2}}{\left({#1}\right)}}

     \providecommand{\constant}[1]{C_{#1}}
     \providecommand{\constext}[1]{\constant{\operatorname{#1}}}            

     \providecommand{\diam}{\operatorname{diam}}
     
     \renewcommand{\div}{\operatorname{div}}

     \providecommand{\interior}{\operatorname{int}}

     \providecommand{\fracl}[2]{{#1}/{#2}}

     \providecommand{\Eye}[1]{
       \begin{bmatrix}
       \ifthenelse{#1>1}{
         \ifthenelse{#1>2}{
           \ifthenelse{#1>3}{
             1&0&\dotso&0
             \\
             0&1&\dotso&0
             \\
             \vdots&\vdots&\ddots&\vdots
             \\
             0&0&\dotso&1
           }{
             1&0&0
             \\
             0&1&0
             \\
             0&0&1
           }
         }{
           1&0
           \\
           0&1
         }
       }{
         1
       }
       \end{bmatrix}
     }
     \providecommand{\Id}{\operatorname{Id}}                   
     
     \providecommand{\maxi}[2]{#1\vee#2}                       
     \providecommand{\maxofset}[1]{\max\setof{#1}}

     \providecommand{\Oh} {\operatorname{O}}                   
     \providecommand{\oh} {\operatorname{o}}                   
     \providecommand{\pd}[2]{\ensuremath{\partial_{#1}{#2}}\xspace} 
     
     \providecommand{\dt}{\ensuremath{\d_t}}
     
     \providecommand{\pdt}[1]{\pd t{#1}}                       
     
     \providecommand{\union}[1]{\ensuremath{\bigcup}_{#1}}

     \renewcommand{\vec}[1]{\ensuremath{\boldsymbol{#1}}}
     \providecommand{\geovec}[1]{\vec{#1}}

      \providecommand{\CC}{\ensuremath{\operatorname C}\xspace}
      \providecommand{\HH}{\ensuremath{\operatorname H}\xspace}
      \providecommand{\LL}{\ensuremath{\operatorname L}\xspace}

      \providecommand{\leb}[1]{\ensuremath{\LL_{#1}}}
     
      \providecommand{\sobh}[1]{\ensuremath{\HH^{#1}}}
      \providecommand{\sobhz}[1]{\sobh{#1}_0}

      \providecommand{\poly}[1]{\ensuremath{\rP}^{#1}}
     
      \providecommand{\fespace}{\rV}
      
      \providecommand{\fes}[1]{\ensuremath{\fespace^{#1}}}
      \providecommand{\fez}[1]{\ensuremath{\fezerospace^{#1}}}

      \providecommand{\EOC}{\ensuremath{\operatorname{EOC}}\xspace}
     
     \providecommand{\Forall}{\:\forall\:}
     
     \providecommand{\Foreach}{\quad\Forall}

     \usepackage{amscd}

     \providecommand{\funk}[3]{\ensuremath{#1:#2\to#3}}

     \renewcommand{\restriction}[2]{\left.#1\right|_{#2}}

     \providecommand{\hhat}[1]{\hat{\hat{#1}}} 
     \providecommand{\ccheck}[1]{\check{\check{#1}}} 
     
     \providecommand{\Program}[1]{\textsf{#1}\xspace}

     \providecommand{\alberta}{\Program{ALBERTA}}

     \providecommand{\secsymbol}{\S}
     
     \providecommand{\secref}[1]{\secsymbol\ref{#1}}

     \providecommand{\eqncomment}[1]{\ensuremath{\qquad\qp{\text{{#1}}}}}

     \providecommand{\ListParameters}{}
     \renewcommand{\ListParameters}
     {
     	 \setlength{\topsep}{0em}
     	 \setlength{\leftmargin}{0em}
              \setlength{\itemsep}{0ex}
     	 \setlength{\parsep}{.5ex}
     	 \setlength{\itemindent}{\labelsep}
     	 \addtolength{\itemindent}{\labelwidth}
     }
     
     \newcounter{LetterListItem}
     \renewcommand{\theLetterListItem}{(\alph{LetterListItem})}
     \newenvironment{LetterList}%
     {
     	\begin{list}%
     	{\theLetterListItem\ }%
     	{\usecounter{LetterListItem}
     	 \ListParameters
     	}
     }%
     {\end{list}}
     \newcounter{NumberListItem}
     \renewcommand{\theNumberListItem}{\arabic{NumberListItem}}
     \newenvironment{NumberList}%
     {
     	\begin{list}%
     	{\theNumberListItem.\ }%
     	{\usecounter{NumberListItem}%
     	 \ListParameters
     	}
     }%
     {\end{list}}
     \newcounter{QuestionListItem}
     \renewcommand{\theQuestionListItem}{\textbf{Question \arabic{QuestionListItem}}}
     {
     	\begin{list}%
     	{\theQuestionListItem.\ }%
     	{\usecounter{QuestionListItem}%
     	 \ListParameters
     	}
     }%
     {\end{list}}
     \newcounter{RomanListItem}
     \renewcommand{\theRomanListItem}{(\roman{RomanListItem})}
     {
     	\begin{list}%
     	{\theRomanListItem\ }%
     	{\usecounter{RomanListItem}
     	 \ListParameters
     	}
     }%
     {\end{list}}
     \newcounter{StepsItem}
     {
     	\begin{list}%
     	{Step \theStepsItem.\ }%
     	{\usecounter{StepsItem}%
     	 \ListParameters
     	}
     }%
     {\end{list}}
     
     \providecommand{\grad}{\nabla}
     \renewcommand{\grad}{\nabla}
     
     \usepackage{amsthm} 
     \ifthenelse{\boolean{isthesis}}{%
       \setcounter{secnumdepth}{1}%
     }{
       \setcounter{secnumdepth}{2} 
     }
     \providecommand{\ListParameters}{}
     \renewcommand{\ListParameters}
     {
     	 \setlength{\topsep}{0em}
     	 \setlength{\leftmargin}{0em}
              \setlength{\itemsep}{0ex}
     	 \setlength{\parsep}{.5ex}
     	 \setlength{\itemindent}{\labelsep}
     	 \addtolength{\itemindent}{\labelwidth}
     }
     \newtheoremstyle{plain}
       {}
       {}
       {\mdseries\slshape}
       {\parindent}
       {\bfseries}
       {.}
       {.5em}
       {}
     
     \newtheoremstyle{note}
       {}
       {}
       {}
       {\parindent}
       {\bfseries}
       {.}
       {.5em}
       {}
     
     \newtheoremstyle{claim}
       {}
       {}
       {\mdseries\slshape}
       {}
       {\bfseries}
       {}
       {.5em}
       {}
     
     \newtheoremstyle{exercise}
       {}
       {}
       {}
       {}
       {\bfseries}
       {.}
       {1em}
       {}
     
     \newtheoremstyle{break}
       {}
       {}
       {}
       {}
       {\bfseries}
       {.}
       {\newline}
       {}
     
       \providecommand{\ObsName}{Remark}
       \providecommand{\RemName}{Remark}
       \providecommand{\NotName}{Notation}
       \providecommand{\BFNName}{Big~Fat~Note}
       \providecommand{\DefName}{Definition}
       \providecommand{\ExaName}{Example}
       \providecommand{\TheName}{Theorem}
       \providecommand{\LemName}{Lemma}
       \providecommand{\ProName}{Proposition}
       \providecommand{\CorName}{Corollary}
       \providecommand{\PbmName}{Problem}
       \providecommand{\HypName}{Hypothesis}
       \providecommand{\AlgName}{Algorithm}
       \providecommand{\ExeName}{Exercise}
       \providecommand{\SolName}{Solution}
       \providecommand{\ClaName}{Claim}
       \providecommand{\EsyName}{Essay}
       \providecommand{\Proofname}{Proof}
     
     \ifthenelse{\boolean{isthesis}}{%
       \providecommand{\Thecounter}{The}
     }{%
       \providecommand{\Thecounter}{subsection}
     }
     \providecommand{\pdfformat}[1]{
        \provideboolean{pdfoutput}
        \setboolean{pdfoutput}{#1}
       \ifthenelse{\boolean{pdfoutput}}{
         \typeout{using pdf}
         \usepackage[pdftex]{graphicx,xcolor}
         \providecommand{\graphext}{pdf}
         \renewcommand{\graphext}{pdf}
         \providecommand{\graphextex}{pdf_t}
         \renewcommand{\graphextex}{pdf_t}
       }{
         \typeout{using eps}
         \usepackage[dvips]{graphicx,xcolor}
         \providecommand{\graphext}{eps}
         \renewcommand{\graphext}{eps}
         \providecommand{\graphextex}{eps_t}
         \renewcommand{\graphextex}{eps_t}
       }
       \usepackage{epsfig}
       \usepackage{tikz}
       \usepackage{rotating}
       \definecolor{SussexFlint}{rgb}{.00,.19,.21}
       \definecolor{SussexGrey}{rgb}{.51,.58,.49}
       \definecolor{SussexOrange}{rgb}{.94,.29,.00}
       \definecolor{SussexYellow}{rgb}{1.00,.73,.00}
       \definecolor{SussexRed}{rgb}{.94,.01,.49}
       \definecolor{SussexPurple}{rgb}{.48,.06,.44}
       \definecolor{SussexGreen}{rgb}{.00,.58,.46}
       \definecolor{SussexBlue}{rgb}{.00,.58,.65}
       \colorlet{a}{SussexOrange}
       \colorlet{b}{SussexRed}
       \colorlet{c}{SussexYellow}
       \colorlet{d}{SussexPurple}
       \colorlet{e}{SussexGreen}
       \colorlet{f}{SussexBlue}
       \colorlet{g}{white}
       \colorlet{h}{SussexGrey}
       \colorlet{i}{black}
       \colorlet{j}{SussexFlint}
       \newcommand{\mausDarkColorTheme}{
         \colorlet{a}{SussexYellow!50!yellow}
         \colorlet{b}{SussexGreen!50!green}
         \colorlet{c}{SussexBlue!50!cyan}
         \colorlet{d}{SussexOrange!50!yellow}
         \colorlet{e}{SussexRed!50!red}
         \colorlet{f}{SussexPurple}
         \colorlet{g}{black}
         \colorlet{h}{SussexFlint!50!black}
         \colorlet{i}{white}
         \colorlet{j}{SussexGrey}
       }
     }
    \RequirePackage{xcolor,graphicx}
    \RequirePackage{ifthen}
    \RequirePackage{epsfig}
    \RequirePackage{tikz}
    \providecommand{\qed}{\vrule height 5pt depth 0pt width 3pt}
    \providecommand{\qqed}{{\raggedright{\ \hfill\qed}}}
    
    \newcounter{passo}

    \newenvironment{Proof}[1][{}]%
    {\par\noindent{\bf \Proofname\ #1}\setcounter{passo}{0}}%
    {\qqed\par}
    \newenvironment{Proof*}[1][{}]%
    {\subsection{\Proofname\ #1}\setcounter{passo}{0}}
    {\qqed\par}
    \swapnumbers{
      \theoremstyle{plain}
      \ifthenelse
          {\boolean{isthesis}}
          {\newtheorem{The}{\TheName}[section]}
          {\newtheorem{The}[subsection]{\TheName}}
     {
       \theoremstyle{plain}

       \newtheorem{Lem}[\Thecounter]{\LemName}
       \newtheorem{Cor}[\Thecounter]{\CorName}
       
       \newtheorem*{The*}{\TheName}
       \newtheorem*{Lem*}{\LemName}
       \newtheorem*{Pro*}{\ProName}
       \newtheorem*{Cor*}{\CorName}
       \newtheorem*{Pbm*}{\PbmName}
       \newtheorem*{Hyp*}{\HypName}
       \newtheorem*{Exe*}{\ExeName}
       \newtheorem*{Txx*}{\ExeName} 
       \newtheorem*{Con*}{Conclusion}
       \newtheorem*{Sum*}{Summary}
     }
     {
       \theoremstyle{claim}

     }
     {
       \theoremstyle{note}
       \newtheorem{Obs}[\Thecounter]{\ObsName}
       
       \newtheorem*{Obs*}{\ObsName}

       \newtheorem{Def}[\Thecounter]{\DefName}
       \newtheorem*{Def*}{\DefName}
       
       \newtheorem*{Exa*}{\ExaName}
     }
    
     {
       \theoremstyle{break}
       
       \newtheorem*{Alg*}{\AlgName}
     }
    }
    
    \swapnumbers
    {
      \theoremstyle{exercise}
        \newcommand{\ExeCounter}{subsection}

      \newtheorem*{Rdn*}{Reading}
      \newtheorem*{Sol*}{\SolName}
    }
    \colorlet{a}{red}
    \colorlet{b}{blue}
    \colorlet{c}{green!50!blue}
    \colorlet{d}{magenta}
    \colorlet{e}{cyan}
    \colorlet{f}{yellow!50!black}
    \colorlet{g}{white}
    \colorlet{h}{black!50}
    \colorlet{i}{black}
    \colorlet{j}{black!75}
\usepackage{natbib}
\usepackage{subfigure}
\setboolean{isamsltex}{true}
\setboolean{shownotes}{true}
\setboolean{showtodo}{true}
         \renewcommand{\Aposteriori}{{A posteriori}\xspace}
         \renewcommand{\aposteriori}{{a posteriori}\xspace}
         
         \newcommand{\apriori}{{a priori}\xspace}
         \newcommand{\estimator}{\cE}
         \newcommand{\ef}[3]{\ensuremath{\estimator[#1,#2,#3]}}
         \newcommand{\inddata}{\ensuremath{\gamma}}
         
         \newcommand{\indspace}{\ensuremath{\eta}}
         \newcommand{\indtime}{\ensuremath{\theta}}
         \newcommand{\indrec}{\ensuremath{\varepsilon}}
         \newcommand{\honew}{{\sobh1(\W)}} 
         \newcommand{\honezw}{{\sobhz1(\W)}} 

         \newcommand{\tn}{{t_n}}
         \newcommand{\tno}{{t_{n-1}}}
         
         \newcommand{\TNO}{{\ensuremath{t_{N-1}}}\xspace}
         \newcommand{\TN}{{\ensuremath{t_N}}\xspace}
         \renewcommand{\TN}{{\ensuremath{T}}\xspace}
         \newcommand{\taun}{\ensuremath{\tau_n}}
         \newcommand{\taunm}{\ensuremath{\tau_n^{-1}}}

         \newcommand{\TAUNO}{\ensuremath{\tau_{N-1}}}
         \newcommand{\TAUN}{\ensuremath{\tau_N}}
         \newcommand{\In}{\ensuremath{(\tno,\tn]}}
         \newcommand{\hn}{\ensuremath{h_n}}
         
         \newcommand{\hathn}{\ensuremath{\hat h_n}}
         
         \newcommand{\fn}{\ensuremath{f^n}}

         \newcommand{\projf}[1]{\ensuremath{{\barf}^{#1}}}

         \newcommand{\U}{\ensuremath{U}}
         
         \newcommand{\un}{\U^n}
         \newcommand{\uno}{\U^{n-1}}
         \newcommand{\UN}{\U^N}
         
         \newcommand{\UZERO}{\U^0}
         
         \newcommand{\wn}{\ensuremath{\w^n}\xspace}
         \newcommand{\wno}{\ensuremath{\w^{n-1}}\xspace}
         
         \newcommand{\contopera}[1]{-\div\left[\vec A\grad #1\right]}
         \newcommand{\contoperax}[2]{-\div\left[\vec A(\vec{#2})\grad #1(\vec#2)\right]}
         \newcommand{\elop}{{\cA}}
         \newcommand{\pwA}{\elop_{\mathrm{el}}}
         \newcommand{\discelop}{A}
         \newcommand{\disca}[1]{\discelop^{#1}}
         \newcommand{\an}{\disca n}
         \newcommand{\ano}{\disca{n-1}}
         \newcommand{\rec}{\cR}

         \newcommand{\lprojop}{P_0}
         \newcommand{\lproj}[1]{\lprojop^{#1}}
         
         \newcommand{\eproj}[1]{P_1^{#1}}
         
         \newcommand{\abil}[2]{\ensuremath{a\left(#1,#2\right)}}

         \newcommand{\bdisc}{\ensuremath{\partial}}

         \newcommand{\fdisc}{\ensuremath{\overset{\smash{\shortrightarrow}}\partial}}
         \newcommand{\sdisc}{\ensuremath{{\partial^2}}}
         \newcommand{\discn}{\ensuremath{\overline\partial}}
         \newcommand{\supfespace}{\ensuremath{\rW}}

\providecommand{\CBS}{the Cauchy--Bunyakovskii--Schwarz inequality\xspace}
\renewcommand{\aposteriori}{aposteriori\xspace}
\renewcommand{\Aposteriori}{Aposteriori\xspace}

\newcommand{\zetas}{\ensuremath{z_s}}
\newcommand{\zetat}{\ensuremath{z_T}}
\renewcommand{\Id}{\operatorname{I}}
\renewcommand{\maxi}[2]{\max\left(#1,#2\right)}

\renewcommand{\fes}[1]{\ensuremath{\tilde\fespace^{#1}}}
\renewcommand{\fez}[1]{\ensuremath{\fespace^{#1}}}
\renewcommand{\rec}[1]{\ensuremath{\cR^{#1}}}
\newcommand{\supcT}{\ensuremath{\cM}}
\setboolean{showchanges}{true}

\newcommand{\f}[1]{\tilde f^{#1}}
\renewcommand{\fn}{{\f n}}
\renewcommand{\projf}[1]{\lproj{#1}\f{#1}}
\renewcommand{\In}{\ensuremath{I_n}}
\newcommand{\indcoarse}[1]{\ensuremath{\gamma_{#1}}}
\renewcommand{\indtime}[1]{\ensuremath{\theta_{#1}}}
\renewcommand{\inddata}[1]{\ensuremath{\beta_{#1}}}
\renewcommand{\indrec}[1]{{\ensuremath{\varepsilon_{#1}}}}
\newcommand{\EI}{\ensuremath{\operatorname{EI}}\xspace}

\numberwithin{equation}{section}
\begin{document}
%
%
\newcommand{\mytitle}%
           {A comparison of duality and energy \aposteriori
             estimates for $\leb\infty(0,T;\leb2(\W))$
             in parabolic problems}
\ifthenelse{\boolean{isamsltex}}%
{ \title[Duality, energy, reconstruction in \aposteriori error control]%
\mytitle \author{Omar Lakkis} 
\address{ Omar Lakkis\newline
Department of Mathematics\newline
University of Sussex\newline
Brighton\newline
GB-BN1 9RF, England UK\newline
\url{http://www.maths.sussex.ac.uk/Staff/OL}
} \curraddr{}
\email{\linkedemail{o.lakkis@sussex.ac.uk}} 
\thanks{This work was partially supported by
the E.U. RTN {\em Hyke} HPRN-CT-2002-00282 and the Marie Curie
Fellowship Foundation. O.L wishes to thank the Hausdorff Institute for Mathematics, Bonn.}
\author{Charalambos Makridakis}
\address{ %
  Charalambos Makridakis\newline
  Department of Applied Mathematics\newline %
  University of Crete\newline %
  GR-71409 Heraklion, Greece\newline
  and
  \newline
  Institute for Applied and Computational Mathematics\newline %
  Foundation for Research and Technology-Hellas\newline %
  Vasilika Vouton P.O.Box 1527\newline %
  GR-71110 Heraklion, Greece}
\curraddr{}
\email{\linkedemail{makr@tem.uoc.gr}}
\author{Tristan Pryer} 
\address{ Tristan Pryer\newline
  School of Mathematics, Statistics \& Actuarial Science\newline 
  University of Kent\newline
  Canterbury\newline 
  GB-CT2 7NF, England UK} 
\curraddr{}
\email{\linkedemail{T.Pryer@kent.ac.uk}} 
\thanks{T.P. was supported at Sussex by a
  EPSRC D.Phil. postgraduate research fellowship.}
\subjclass[2000]{Primary: 65N30}
\date{\today}
\commby{}
\dedicatory{}
}
{
\title{\mytitle}
\author{Omar Lakkis, Charalambos Makridakis and Tristan Pryer}
\date{\today}
}
\maketitle
\begin{abstract}
We use the elliptic reconstruction technique in combination with a
duality approach to prove \aposteriori error estimates for fully
discrete backward Euler scheme for linear parabolic equations.  As an
application, we combine our result with the residual based estimators
from the \aposteriori estimation for elliptic problems to derive
space-error indicators and thus a fully practical version of the
estimators bounding the error in the $\leb\infty(0,T;\leb2(\W))$
norm. These estimators, which are of optimal order, extend those
introduced by \citet{ErikssonJohnson:91} by taking into account the
error induced by the mesh changes and allowing for a more flexible use
of the elliptic estimators.  For comparison with previous results we
derive also an energy-based \aposteriori estimate for the
$\leb\infty(0,T;\leb2(\W))$-error which simplifies a previous one
given in \cite{lakkis-makridakis:06}.  We then compare both estimators
(duality vs. energy) in practical situations and draw conclusions.
\end{abstract}
%
%
\section{Introduction}
\Aposteriori error estimators and their use to derive adaptive mesh
refinement algorithms to solve time-dependent problems constitute the
object of current research.  The problem is appealing for the
theoretician as a test ground for novel analytical techniques as well
as for the practitioners which are interested in minimizing the amount
of computational time in order to obtain a satisfactory accuracy in
the computer simulations of time-dependent PDE's.  Both the
theoretical and practical aspects of \aposteriori-based adaptive
numerical methods for evolution partial differential equations has
benefited immensely from the surge in the production of dedicated
papers in the last 20 years, although fundamental questions such as
convergence of adaptive algorithm remains open.

In this paper, we address the problem of \aposteriori error
estimation for the time-dependent model problem
\begin{equation}
  \pdt u(\vec x,t)+\cA u(\vec x,t)=f(\vec x,t)
\end{equation}
for $\vec x\in\W\subseteq\R d$ and $0\leq t\leq T$, where $\cA$ is an
elliptic operator to be described in detail further in
\S\ref{sec:setup}.  In a previous article
\citep{lakkis-makridakis:06}, we used the \emph{elliptic
  reconstruction} in combination with \emph{energy techniques} to
derive \aposteriori error estimates the heat equation in
$\leb\infty(0,T;\leb2(\W))$'s norm, to analyze fully discrete implicit
Euler method in time and conforming finite element methods (FEM) in
space.  The elliptic reconstruction was then used alongside a
parabolic energy technique to derive optimal-order \aposteriori
residual-based $\leb\infty(0,T;\leb2(\W))$-error estimators.  Previous
work for the spatially semidiscrete scheme was introduced by
\citet{makridakis-nochetto:03}.  The elliptic reconstruction has since
then been used later as an analytical tool, in combination with energy
or other techniques to deal with time, in order to establish estimates
in various norms for linear and nonlinear problems
\citep[e.g.]{ bartels-muller:09:preprint, demlow-lakkis-makridakis:09,
  DemlowMakridakis:10, ErnMeunier:09, GeorgoulisLakkis:10}.  

Before the introduction of the elliptic reconstruction,
$\leb\infty(0,T;\leb2(\W))$-error estimates could be derives by using
the \emph{duality techinque}, at the cost of assuming restrictive
assumptions on the domain (e.g., convexity) as well as on the mesh
\citep{ErikssonJohnson:91}.
Our chief goals in this paper are
\begin{NumberList}
\item 
  to explore, for the first time, the possibility of using the
  \emph{elliptic reconstruction} technique in conjunction with the
  \emph{duality} technique as introduced by
  \citet{ErikssonJohnson:91};
  \\
  and
\item
  to compare duality estimates with \emph{energy estimates} for the
  same norm; here the use of the elliptic reconstruction is crucial as
  it provides a simple abstract result for the $\leb\infty(\leb2)$ norm
  which is the same that is used in duality.
\end{NumberList}

The \emph{duality} technique provides an important alternative to
energy techniques and is widely used for the derivation of \apriori
and \aposteriori error estimates both for elliptic and parabolic
problems.  Since being first considered by it has been developed in
many different directions, including its use in \emph{implicit and
  goal oriented \aposteriori error estimates}.

The elliptic reconstruction has been used in combination with energy
estimates, where one mimics the energy estimates for the parabolic
equation in order to derive error estimates from a PDE where the
error, or part thereof, is the ``unknown''.  In this paper, we exhibit
the flexibility of the elliptic reconstruction technique by showing
that it can be completely decoupled from energy considerations (or any
other method used to deal with time integration and time-stepping, for
that matter).  This is not obvious, indeed, in many works \aposteriori
analysis, the elliptic part is entangled with the parabolic part and
there is not a clear cut difference between elliptic and parabolic
effects.  As noted in recent work on \aposteriori analysis for
time-dependent problems \citep[e.g.]{ AkrivisMakridakisNochetto:06,
  BergamBernardiMghazli:05, BernardiVerfurth:04, FrutosNovo:02,
  picasso:98} understanding the splitting between the elliptic,
stationary, and parabolic, time-dependent, errors, as well as the part
of the error where these effects are coupled, is important in
designing adaptive methods and avoiding repetition.

An important by-product of our approach is that the mesh-change in
time is considered as part of the proofs of our theorems. Indeed,
unlike former derivations \aposteriori error estimates via duality
\citep[mainly]{ErikssonJohnson:91}, we do not impose on the mesh any
assumption that are susceptible of violation in a practical
implementation of the scheme, such as the no-refinement assumptions.

From a more practical side, we give an application of our theory, by
comparing in a series of benchmarks where elliptic $\leb2(\W)$
residual-based estimators are used
\citep{ainsworth-oden:book,liao-nochetto:03}.  We emphasize, however,
that our results are not limited to the use of residual-based
estimators and that other estimators which work for the $\leb2(\W)$ norms
in elliptic problems could be used
\citep[e.g.]{LakkisPryer:10}.

Our main results in this paper are duality-based estimates,
Theorem~\ref{the:general.aposteriori.estimate} and
Corollary~\ref{cor:fully-discrete-estimate}, an energy-based esimate,
Theorem~\ref{the:general.energy.aposteriori.estimate} and a computer
experiment designed at comparing in practice both estimators.
From a theoretical perspective, Corollary
\ref{cor:fully-discrete-estimate} generalizes the duality estimates of
\citet{ErikssonJohnson:91}, mainly by removing unrealistic assumptions
on the meshes.  A direct application of the duality-estimate Theorem
\ref{the:general.aposteriori.estimate} provides finer estimates with
respect to time accumulation. This is especially helpful in situations
where the error (on a time-invariant mesh) decreases with time and for
long-time integration.  Finally, energy-estimate
Theorem~\ref{the:general.energy.aposteriori.estimate} simplifies (by
using PoincarÃ© inequality) special cases from
\citet{lakkis-makridakis:06} and provides the basis for our
comparison.
The numerical results show that the estimators behave roughly the
same, with a slight edge for the energy-based ones when it comes to
time accumulation and long time integration.  This is a confirmation
of the theoretical observation that the ``tails'' of the coefficients
for the time-accumulation are much heavier for the duality estimators
(see Figure \ref{fig:time-accumulation-coefficients}).  The energy
estimator benefits from an exponential decay in these coefficients
which also provides a faster way of computing then and a more
economical storage.  In summary, we found that if energy estimators
are available they are better suited for practical scenarios where the
$\leb\infty(0,T;\leb2(\W))$ is important.

The rest of this article is organized as follows: In \S\ref{sec:setup}
we recall the main tools related to the elliptic reconstruction. In
\S\ref{sec:duality.semi} we analyze the spatially semidiscrete scheme
using a duality approach.  In \S\ref{sec:duality.fully} we extend the
\S\ref{sec:duality.semi} to the fully discrete scheme and in
\S\ref{sec:proof-duality} give the proof of those results. In
\S\ref{sec:energy} we state and prove the estimates based on the
energy approach.  Finally in \S\ref{sec:sample-application-residuals}
we summarize our computer experiments from which we drew the main
practical conclusions of this research.
\subsection*{Acknowledgments}
Part of this research is based on work from O.L. stay at FORTH in
Crete in the framework of a Marie Curie fellowship and the HYKE
RTN. T.P. was funded mostly by an EPSRC postgraduate research
fellowship during this research.  O.L. thanks Christoph Ortner and
SÃ¶ren Bartels for the interesting discussions about elliptic
reconstruction and the energy approach at the Hausdorff Institute for
Mathematics, Bonn, which we also thank for its kind generosity and
outstanding hospitality.
\section{The discrete scheme and the elliptic reconstruction}
\label{sec:setup}

In this section we introduce the numerical schemes that we study,
some basic tools including the definition of the elliptic
reconstruction.

\subsection{Basic set-up}
\label{sse:setting.notation}
We introduce next the PDE whose discretization is the object of this
paper.  Let $\W$ be a bounded domain of the Euclidean space $\rR^d$,
for some fixed positive integer \emph{space dimension} $d$ and a
\emph{final time} $T\in\rR^+$.  We shall assume throughout this
paper's discussion that \W is a polygonal convex domain, noticing that
all the results can be extended to certain non-convex domains, like
domains with reentrant corners in $d=2$, following ideas of
\citet{liao-nochetto:03} regarding the elliptic \aposteriori
$\leb2(\W)$-error estimates.

Given a Lebesgue measurable set $D\subset\rR^d$, we define
\begin{align}
  \ltwop\phi\psi_D
  &
  :=
  \int_D\phi(\vec x)\psi(\vec x)\mu(\!\d\vec x),
  \\
  \Norm{\phi}_D
  &
  :=
  \Norm{\phi}_{\leb2(D)}
  :=
  \ltwop\phi\phi_D^{1/2},
  \\
  \norm{\phi}_{k,D}
  &
  :=
  \Norm{\D^k\phi}_D,\text{ for }k\in\Nat
  \\
  \Norm{\phi}_{k,D}
  &
  :=
  \bigg(\Norm{\phi}_D^2+\sum_{j=1}^k\norm{\phi}_{j,D}^2\bigg)^{1/2},
  \text{ for }k\in\Nat,
\end{align}
where $\mu(\!\d\vec x)$ denotes either the Lebesgue measure element $\d\vec
x$, when $D$'s such measure is positive, or the $(d-1)$-dimensional
(Hausdorff) measure $\ds{\vec x}$, when $D$ has zero Lebesgue measure.
In many instances, in order to compress notation and when there is no
danger of engendering confusion, we may drop altogether the
``differential'' symbol from integrals. This convention applies also
to integrals in time.

We will use the standard~\citep{evans:PDE} function spaces $\leb2(D)$,
$\sobh k(D)$, $\sobhz k(D)$ and denote by $\sobh{-1}(D)$ the dual
space of $\sobhz1(D)$ with the corresponding pairing written as
$\duality{\cdot}{\cdot}_D$.  We omit the subscript $D$ whenever
$D=\W$.  We denote the PoincarÃ©--Friedrichs constant associated with $\W$
by $\constext{PF}$ and \emph{we take the seminorm $\norm{\cdot}_1$ to
  be the norm of $\honezw$}.  We use the usual duality identification
  \begin{equation}
    \honezw\subset\leb2(\W)\sim\leb2(\W)'\subset\sobh{-1}(\W)
  \end{equation}
and the dual norm
\begin{equation}\label{eqn:dual.norm}
  \Norm\psi_{-1}
  :=\sup_{0\neq\phi\in\sobhz1(\W)}\frac{\duality\psi\phi}{\norm{\phi}_1}
  \left(=\sup_{0\neq\phi\in\sobhz1(\W)}\frac{\ltwop\psi\phi}{\norm{\phi}_1},
  \text{ if }\psi\in\leb2(\W)\right).
\end{equation}

Let $a$ be the elliptic bilinear form defined on $\honezw$ by
\begin{equation}
  \label{eqn:bilinear.form}
  \abil v\psi :=\ltwop{\vec A\grad v}{\grad\psi}\Foreach v,\psi\in\honezw,
\end{equation}
where ``$\grad$'' denotes the spatial gradient and the matrix-valued
function $\vec A\in\leb\infty(\W)^{d\times d}$ is such that
\begin{gather}
  \label{eqn:bounded.bilinear}
  \abil\psi\phi \leq \beta\norm{\psi}_1\norm{\phi}_1
  \Foreach\phi,\psi\in\honezw,\\
  \label{eqn:coercive.bilinear}
  \abil\phi\phi \geq \alpha\norm{\phi}_1^2
  \Foreach\phi\in\honezw,
\end{gather}
with $\alpha,\beta\in\reals^+$.
We also use the {\em energy norm} $\norm{\cdot}_a$
defined as
\begin{equation}
  \norm{\phi}_a:=\abil\phi\phi^{1/2}\Foreach\phi\in\honezw.
\end{equation}
It is equivalent to the norm $\norm{\cdot}_1$ on the space $\sobhz1(\W)$, in
view of \eqref{eqn:bounded.bilinear} and
\eqref{eqn:coercive.bilinear}.  In particular, we will often use the following
inequality
\begin{equation}
  \label{eqn:one.energy.inequality}
  \norm\phi_1\leq\alpha^{-1/2}\norm\phi_a\Foreach\phi\in\honezw.
\end{equation}

Let $u\in\leb\infty(0,T;\honezw)$, with
$\partial_tu\in\leb2(0,T;\sobh{-1}(\W))$, be the unique solution of the
linear parabolic problem
\begin{equation}
  \begin{split}\label{eqn:continuous.heat}
    &\duality{\partial_t{u}}{\phi}+\abil u\phi=\ltwop{f}{\phi}
    \Foreach\phi\in\honezw,\\
    &\text{ and }u(0)=g,
  \end{split}
\end{equation}
where $f\in\leb2(\W\times(0,T))$ and $g\in\honezw$.   
Whenever not stated explicitly, we assume that the data $f,g,\vec A$
and the solution $u$ of the above problem are sufficiently regular for
all the norms involved to make sense.

In order to discretize the time variable in
\eqref{eqn:continuous.heat}, we introduce the partition
$0=t_0<t_1<\dotsb<t_N=T$ of $[0,T]$.  Let $I_n:=\opclinter\tno\tn$ and we
  denote by $\tau_n:=t_n-t_{n-1}$ the time steps. We will consistently
  use the following ``superscript convention'': whenever a function
  depends on time, e.g. $f(\vec x, t)$, and the time is fixed to be
  $t=t_n,\,\rangefromto n0N$ we denote it by $f^n(\vec x)$.  Moreover,
  we often drop the space dependence explicitly, e.g, we write $f(t)$
  and $f^n$ in reference to the previous sentence.

We use a conforming fixed polynomial degree FEM to discretize the
space variable.  Let $(\cT_n)_{\rangefromto n0N}$ be a family of
conforming triangulations of the domain $\W$
\citep{brenner-scott,ciarlet:fem}.  These triangulations are allowed to
change at each timestep, as long as they stay {\em
  compatible}~\citep[\S A]{lakkis-makridakis:06}, which is an extremely
mild requirement automatically implemented by many refinement methods.
%

For each given a triangulation $\cT_n$, we denote by $h_n$ its
meshsize function defined as
\begin{equation}
  h_n(\vec x)=\diam(K),\text{ where }{K\in\cT_n}\text{ and }{\vec x\in K},  
\end{equation}
for all $\vec x\in\W$.
We also denote by $\cS_n$
the set of {\em internal sides} of $\cT_n$, these are edges in
$d=2$---or faces in $d=3$---that are contained in the interior of
\W; the {\em interior mesh of edges} $\Sigma_n$ is then defined as
the union of all internal sides $\cup_{E\in\cS_n}E$.  We associate
with these triangulations the {\em finite element spaces}:
\begin{gather}
  \label{eqn:def:finite-element-spaces}
  \fez n:=\ensemble{\phi\in   \honezw}
       {\forall K\in\cT_n:\restriction\phi K\in\poly\ell},
\end{gather}
where $\poly\ell$ is the space of polynomials in $d$ variables of
degree at most $\ell\in\Nat$.
Given two successive compatible triangulations $\cT_{n-1}$ and
$\cT_n$, we define $\hat h_n:=\maxi{h_n}{h_{n-1}}$
\citep[Appendix]{lakkis-makridakis:06}.  
We will also use the sets
$\hat\Sigma_n:=\Sigma_n\cap\Sigma_{n-1}$ and
$\check\Sigma_n:=\Sigma_n\cup\Sigma_{n-1}$.  
To keep notation light, we shall often use two ``generic'' finite
element spaces $\fespace$ and $\supfespace$, defined as in
(\ref{eqn:def:finite-element-spaces}) in association with two
``generic'' triangulations $\cT$ and $\supcT$, respectively.

\begin{Def}[fully discrete scheme]
  \label{def:fully.discrete}
  We consider the following fully discrete scheme of problem
  \eqref{eqn:continuous.heat} associated with the finite element
  spaces $\fez n$:
  \begin{equation}
    \label{eqn:fully.discrete.scheme}
    \begin{gathered}
      U^0:=I^0u(0),\text{ and }
      \\
      \taunm\ltwop{\un-\uno}{\Phi_n}
      +\abil{\un}{\Phi_n}=\ltwop{\f n}{\Phi_n}
      \Foreach\Phi_n\in\fez n,\text{for $\rangefromto n1N$}.
    \end{gathered}
  \end{equation}
  Here the operator $I^0$ is some suitable interpolation or projection
  operator from $\honezw$, or $\leb2(\W)$, onto $\fez n$, and $\f n$
  equals either the value of $f$ at $\tn$, $f^n:=f(\cdot,\tn)$, or its
  time-average on $\In$, $\int_\tno^\tn f(\cdot,t)\d t/(\tn-\tno)$.
  This scheme is the standard backward (or implicit) Euler--Galerkin
  finite element scheme~\citep{Thomee:2006:book}.
\end{Def}
%
%
In the sequel we shall use a continuous piecewise linear extension in
time of the sequence $(\un)$ which we denote by $U(t)$ for $t\in[0,T]$
(see \secref{def:discrete.time} for the precise definition).

\subsection{\Aposteriori estimates and reconstruction operators.}
The elliptic reconstruction, as described by Makridakis \&
Nochetto~\citep{makridakis-nochetto:03} consists in associating with
$\funk U{\clinter0T}{\fespace}$ an auxiliary function
$\funk\w{\clinter0T}\honezw$, in such a way that when the {\em total error}
\begin{equation}
  e:=U-u
\end{equation}
is decomposed as follows
\begin{gather}
  e=\rho-\epsi\\
  \epsi:=\w-U,\quad\rho:=\w-u,
\end{gather}
then the following properties are satisfied:
\begin{NumberList}
\item The error $\epsi$ is easily controlled by elliptic \aposteriori quantities of optimal order.
\item The error $\rho$ satisfies  a modification of the original PDE whose right-hand side depends on 
$\epsi$ and $U$.  This right-hand side can be bounded \aposteriori in an optimal way.
\end{NumberList}
Therefore in order to successfully apply this idea we must select a
suitable reconstructed function $\w$.  In our case, this choice is
dictated by the elliptic operator at hand; the precise definition is
given in \secref{def:elliptic.reconstruction}.  In addition the effect
of mesh modification will reflect in the right-hand side of the
equation for $\rho$.  As a result of our choice for $\w$ we are able
to derive optimal order estimators for the error in
$\leb\infty(0,T;\leb2(\W))$, as well as in $\leb\infty(0,T;\honezw)$ and
$\sobh1(0,T;\leb2(\W))$. In addition, our choosing $\w$ as the elliptic
reconstruction will have the effect of separating the spatial
approximation error from the time approximation as much as possible.
We show that the spatial approximation is embodied in $\epsi$ which
will be referred to as the {\em elliptic reconstruction error} whereas
the time approximation error information is conveyed by $\rho$, a fact
that motivates the name {\em main parabolic error} for this term.
This ``splitting'' of the error is already apparent in the spatially
discrete case \citep{makridakis-nochetto:03}.

With the above notation, we prove in the sequel that $\rho$ satisfies the following 
variational equation.
\begin{Lem}[main parabolic error equation]
\label{lem:main.error.equation}
For each $\rangefromto n1N$, and for each $\phi\in\sobhz1(\W)$,
\begin{equation}\label{eqn:main.error.equation}
\begin{split}
  \duality{\partial_t\rho}{\phi}+\abil{\rho}{\phi}
  =&\ltwop{\partial_t\epsi}{\phi}+\abil{\w-\w^n}{\phi}\\
   &+\ltwop{\lproj n{f}^n-f}{\phi}
  +\taunm\ltwop{\lproj n\uno-\uno}{\phi}
  \text{ on $I_n$.}
\end{split}
\end{equation}
Here $\lproj{n}$ denotes the $L^{2}$-projection into  $\fez n .$
\end{Lem}

Since the definitions some in parts of this section are independent of
the time discretization and could be applied to any finite element
space, in this section we use two generic $\honezw$-conforming
Lagrange finite element spaces $\fespace$ and $\supfespace$.

Whenever $\fespace$, or $\supfespace$, coincides with one of the $\fez
n$ introduced in , we replace all indexes $\fespace$ by $n$.
\begin{Def}[representation of the elliptic operator, 
    discrete elliptic operator, projections]
  \label{obs:representation}
  Suppose a function $W\in\supfespace$, the bilinear form can be
  then represented as
  \begin{equation}\label{eqn:representation.long.notation}
    \abil v\phi = 
    \sum_{K\in\cT}
    \ltwop{\contopera{W} }{\phi}_K +
    \sum_{E\in\cS}\ltwop{J[W]}\phi_E
    \Foreach\phi\in\honezw,
  \end{equation}
  where $J[W]$ is the {\em spatial jump of the field
 $\vec A\grad W$ across
  an element side $E\in\cS$} defined as
  \begin{equation}
    \begin{split}
      \restriction{J[W]}{E}(\vec x)
      =&\jump{\vec A\grad W}_E(\vec x)
      \\
      :=&\lim_{\ep\rightarrow0}
	\qb{
	  \vec A(\vec x)\grad W(\vec x+\ep\vec\nu_E(\vec x))
	  -\vec A(\vec x)\grad W(\vec x-\ep\vec\nu_E(\vec x))
	}
	\inner\vec\nu_E(\vec x)
    \end{split}
  \end{equation}
  where $\vec\nu_E$ is a choice, which does not influence this
  definition, between the two possible normal vectors to $E$ at the
  point $\vec x$.

  Since we use the representation
  \eqref{eqn:representation.long.notation} quite often, we introduce
  now a practical notation that makes it shorter and thus easier to
  manipulate in convoluted computations.  For a finite element function,
  $W\in\supfespace$ (or more generally for any Lipschitz continuous
  function $w$ that is $\CC^2(\interior(K))$, for each $K\in\cT$), denote by
  $\pwA W$ the {\em regular part} of the distribution $\contopera{{W}}$,
  which is defined as a piecewise continuous function such that
  \begin{equation}
    \ltwop{\pwA W}\phi = \sum_{K\in\cT}\int_K{\contoperax{W} x}\phi(\vec x)\d\vec x
    \Foreach\phi\in\honezw.
  \end{equation}
  The operator $\pwA{} $ is sometime referred to, in the finite
  element community, as the \emph{elementwise elliptic operator}, as
  it can be viewed as the result of the application of
  $\contopera\cdot$ only on the interior of each element $K\in\cT$.
  Although this is a misnomer (as the operator itself does not depend
  in any way on the finite element space) this observation justifies
  our subscript in the notation.  We shall write the representation
  \eqref{eqn:representation.long.notation} in the shorter form
    \begin{equation}\label{eqn:representation.short.notation}
      \abil {W}\phi = \ltwop {\pwA {W}}\phi + \ltwop{J[W]}\phi_{\Sigma}
      \Foreach\phi\in\honezw,
    \end{equation}
    where $\Sigma=\union{E\in\cS}E$.
\end{Def}

Let us now recall some more basic definitions that we will be using.
The {\em discrete elliptic operator} associated with the bilinear
form $a$ and the finite element space $\fespace$ is the operator
$\disca{\fespace}:\honezw\rightarrow\fespace$ defined by
\begin{equation}
  \label{eqn:discrete.elliptic.operator}
  \ltwop{\disca{\fespace}v}{\Phi}=\abil{v}{\Phi}\Foreach\Phi\in\fespace,
\end{equation}
for $v\in\honezw$.

The {\em $\leb2$-projection operator} is defined as the operator
$\lproj{\fespace}: \leb2(\W)\rightarrow\fespace$ such that 
\begin{equation}
  \ltwop{\lproj{\fespace}{v}}{\Phi}=\ltwop{v}{\Phi}
  \Foreach\Phi\in\fespace,
\end{equation} 
for $v\in\leb2(\W)$; and the {\em elliptic projection operator}
$\eproj{\fespace}: \honezw\rightarrow\fespace$ is defined by
\begin{equation} 
  \abil{\eproj{\fespace}{v}}{\Phi}=\abil{v}{\Phi}
  \Foreach\Phi\in\fespace.
\end{equation}
\begin{Def}[elliptic reconstruction]
  \label{def:elliptic.reconstruction}
  We define the {\em elliptic reconstruction operator} associated with
  the bilinear form $a$ and a given finite element space $\fespace$ to be
  the unique operator $\rec{\fespace}:\honezw\rightarrow\honezw$ such that
  \begin{equation}\label{eqn:reconstruction}
    \abil{\rec{\fespace}{v}}{\phi}=\ltwop{\disca\fespace v}{\phi}\Foreach\phi\in\honezw,
  \end{equation} for each given $v\in\honezw$. The function $\rec\fespace v$ is
  referred to as the {\em elliptic reconstruction} of $v$.  

  Note that the domain of the reconstruction operator $\rec\fespace$
  can be taken to be $\honezw$, but it will be used effectively on the
  finite element space and we generally consider its restriction to
  $\fespace$.  The elliptic reconstruction operator $\rec\fespace$,
  restricted to the space $\fespace$, is a right, but not left,
  inverse of the well-known~\cite{Thomee:2006:book} elliptic (or Ritz) projection.
\end{Def}
\begin{Obs}[Galerkin orthogonality]
  A crucial property of the elliptic reconstruction operator $\rec n$
  is that for $v\in\honezw$, $v-\rec n{v}$ is
  $\abil\cdot\cdot$-orthogonal to $\fez n$, i.e.,
  \begin{equation}\label{eqn:reconstruction.orthogonality}
    \abil{v-\rec n v}{\Phi}=0\Foreach\Phi\in\fez n.
  \end{equation}
  This is known as the \emph{Galerkin orthogonality} of the error in
  the finite element literature and is the crucial property that allows
  to obtain \apriori and \aposteriori error estimates.
\end{Obs}
\begin{Def}[elliptic \aposteriori error estimator functional]
  \label{def:elliptic-error-estimator-functional}
  Given a normed functional space $\cV$ containing $\honezw$, (e.g.,
  $\cV=\leb2(\W)$ or $\honezw$) and a generic finite dimensional
  subspace $\rV$, we call \emph{estimator functional} associated with
  the bilinear form $a$, defined in (\ref{eqn:bilinear.form}), the
  space $\rV$ in the the norm of $\cV$, a functional of the form
  \begin{equation}
    \funk{\ef\cdot\rV\cV}{\rV}\reals
  \end{equation}
  such that for each $V\in\rV$ we have
  \begin{equation}
    \Norm{V-\rec\fespace}_{\cV}\leq\ef V\rV\cV.
  \end{equation}
  Thanks to many different techniques
  \cite{ainsworth-oden:book,braess:book,verfuerth:book}, it is
  well-known that there exist many such functionals.  One of the
  simplest examples is given by the residual-based estimator
  functional, justified next by
  Lemma~\ref{lem:elliptic.aposteriori.estimates}, which we will use in
  this work, but we note that our approach can be easily adapted to
  accommodate other estimators.
\end{Def}
\begin{Lem}[residual-based \aposteriori error estimates]
  \label{lem:elliptic.aposteriori.estimates}
  Let $\fespace$ be a finite element space on a triangulation $\cT$
  with edge set $\Sigma$ of the polygonal domain $\W$ as defined in
  \S\ref{sse:setting.notation}.  For any $V\in\fespace$ we have
  \begin{align}
    \label{eqn:elliptic.estimate.h1}
    \norm{\rec\fespace{V}-V}_1
    &
    \leq\frac{C_{3,1}}\alpha
    \Norm{(\pwA V-\an V)\hn}+\frac{C_{5,1}}\alpha
    \Norm{J[V]\hn^{1/2}}_{\Sigma_n},
    \intertext{ and, if furthermore $\W$ is convex, then }
    \label{eqn:elliptic.estimate.l2}
    \Norm{\rec\fespace{V}-V}
    &
    \leq C_{6,2}\Norm{(\pwA V-\an V)\hn^2}+C_{10,2}\Norm{J[V]\hn^{3/2}}_{\Sigma_n},
  \end{align}
  for the $\alpha$ given by (\ref{eqn:one.energy.inequality}) and some
  ($V$-independent) constants $C_{k,j}$, defined in \cite[Appendix~B]{lakkis-makridakis:06}.
\end{Lem}
\begin{Def}[discrete time extensions and derivatives]
\label{def:discrete.time}
Given any discrete function of time---that is, a sequence of
values associated with each time node $t_n$---e.g., $(\un)$, we
associate to it the continuous function of time defined by the
Lipschitz continuous piecewise linear interpolation, e.g.,
\begin{equation}
  \label{eqn:pwl.extension}
  U(t):=l_{n-1}(t)\uno+l_n(t)\un,
  \text{ for $t\in I_n$ and $\rangefromto n1N$};
\end{equation}
where the functions $l_n$ are the hat (linear Lagrange basis) functions defined by
\begin{equation}
  \label{eqn:def:time-lagrange-basis}
  l_n(t):=\frac{t-t_{n-1}}{\tau_n}\one_{I_n}(t)-\frac{t-t_{n+1}}{\tau_{n+1}}\one_{I_{n+1}}(t),
  \text{ for $t\in[0,T]$ and $\rangefromto n0N$},
\end{equation}
$\one_X$ denoting the characteristic function of the set $X$.  

In the sequel will use the following shorthand
\begin{equation}
  \w^\fespace = \rec\fespace U\eqncomment{and thus $\wn=\rec n\un$},
\end{equation}
to denote the elliptic reconstruction of the (semi-)discrete
solution $U$ and $\un$.

The {\em time-dependent elliptic reconstruction} of $U$ is the function
\begin{equation}
  \w(t):=l_{n-1}(t)\rec {n-1}\uno+l_n(t)\rec n\un,
  \text{ for $t\in I_n$ and $\rangefromto n1N$},
\end{equation}
which results in a Lipschitz continuous function of
time.

We introduce next time-discrete derivative (i.e., difference) operators:
\begin{LetterList}
\item 
  {\em Discrete (backward) time derivative}
  \begin{equation}
    \bdisc\un:=\frac{\un-\uno}{\taun}.
  \end{equation}
Notice that $\bdisc\un=\partial_t U(t)$, for all $t\in I_n$, hence we
can think of $\bdisc\un$ as being the value of a discrete function at
$t_n$.  We thus define $\bdisc U$ as the piecewise linear extension of
$(\bdisc\un)_n$, as we did with $U$.
\item {\em Discrete (centered) second time derivative}
  \begin{equation}
    \sdisc\un:=\frac{\bdisc U^{n+1}-\bdisc\un}{\taun}.
  \end{equation}
\item  
  \label{itm:averaged-time-derivative}
  {\em Averaged ($\leb2$-projected) discrete time derivative}
  \begin{equation}
    \label{eqn:averaged.discrete.time.derivative}
    \discn\un := \lproj n\bdisc\un = \frac{\un-\lproj n\uno}{\taun}
    \Foreach
    \rangefromto n 1N.
  \end{equation}
  This last definition stems from $\bdisc\un$ not necessarily
  belonging to $\fez n$ (e.g., when $\fez{n-1}\not\subseteq\fez n$ ), whereas
  $\discn\un\in\fez n$ is always satisfied.
\end{LetterList}
\end{Def}
\begin{Obs}[pointwise form]
  The discrete elliptic operators $\an$ can be employed to write the
  fully discrete scheme \eqref{eqn:fully.discrete.scheme} in the
  following {\em pointwise form}
  \begin{equation}
    \label{eqn:pointwise.form}
    \discn\un(\vec x) + \an\un(\vec x) = \projf n(\vec x)\Foreach\vec x\in\W.
  \end{equation}
  Indeed, in view of $\discn\un+\an\un-\projf n\in\fez n$,
  \eqref{eqn:fully.discrete.scheme}, and
  \eqref{eqn:discrete.elliptic.operator}, we have
  \begin{equation}
    \begin{split}
    \ltwop{\an\un+\discn\un-\projf n}{\phi}
    &=\ltwop{\an\un+\discn\un-\projf n}{\lproj n\phi}\\
    &=\abil{\un}{\lproj n\phi}+\ltwop{\taunm(\un-\uno)-f^n}{\lproj n\phi}
    =0,
  \end{split}
\end{equation}
for any $\phi\in\honezw$. 
Therefore the function $\discn\un+\an\un-\projf n$ vanishes.
\end{Obs}

\subsection{Error equation}
\label{sec:error-equation}
Let us consider the {\em (full) error}, the {\em elliptic
  reconstruction error} and the {\em parabolic error} which are
defined, respectively as follows
\begin{gather}
  e = U - u,\\
  \epsi = \w - U,\\
  \rho = \w - u.
\end{gather}
We have the following decomposition of the error
\begin{equation}
  e = \rho-\epsi.
\end{equation}
We can also readily derive the following error relation for the
parabolic error in terms of the reconstruction error and the
reconstruction itself \cite{lakkis-makridakis:06}:
\begin{equation}
  \label{eqn:error.equation}
  \begin{split}
    \ltwop{\partial_t\rho(t)}\phi +\abil{\rho(t)}\phi
    =&
    \ltwop{\partial_t\epsi(t)}\phi
    +\abil{\w(t)-\wn}\phi\\
    &+\taunm\ltwop{\lproj n\uno-\uno}\phi
    +\ltwop{\projf n-f(t)}\phi
  \end{split}
\end{equation}
for all $\phi\in\honezw$, $t\in I_n$ and $\rangefromto n1N$.
\section{A duality--reconstructive derivation of \aposteriori error estimates}
\label{sec:duality.semi}
In this section we synthetically describe how the combination of the
elliptic reconstruction and the parabolic duality techniques provides
\aposteriori error estimates.  To keep the discussion as simple as
possible, we study first the \emph{spatially semidiscrete scheme}.
This simplification allows us to expose our main ideas, which we
employ later for the fully discrete case in
\S\ref{sec:duality.fully}.
\subsection{Notational warning}
\label{sec:notational-warning}
Since we will be dealing with the space semidiscrete scheme only, we
will use the same symbols introduced for the fully discrete scheme in
\secref{sec:setup}, albeit in their semidiscrete analog by dropping
the index $n$.  \emph{The notation now introduced is valid only in
  this section.}  In particular time-dependent functions, such as $U$,
$\w$, $e$, $\epsi$ and $\rho$, to be introduced next, should not be
confused with their fully-discrete analogs introduced earlier in
\S\ref{sec:setup} and valid outside this section.
\subsection{Notation, spatially semidiscrete scheme and the error relation}
\label{sec:notation-and-semidiscrete-error-relation}
Let $\fespace$ be a given (time-invariant) finite element space, as
defined in \S\ref{sse:setting.notation}, consider the function $\funk
U{\clinter0T}{\fespace}$ which satisfies the following
\emph{semidiscrete Galerkin finite element scheme} associated with
the PDE (\ref{eqn:continuous.heat}):
\begin{equation}
  \label{eqn:semidiscrete.scheme}
  \begin{gathered}
    U(0):=I u(0),\text{ and }
    \\
    \ltwop{\partial_t U(t)}{\Phi}+\abil{U(t)}{\Phi}
    =
    \ltwop{f(t)}\Phi\Foreach\Phi\in\fespace,t\in[0,T],
  \end{gathered}
\end{equation}
where the operator $I$ is a suitable interpolation or projection
operator from $\honew$, or $\leb2(\W)$, onto $\fespace$.  

We define the (full) \emph{error} at time $t$ to be $e(t):=U(t)-u(t)$
and the \emph{semidiscrete} elliptic reconstruction to be
$\w(t):=\rec\fespace U(t)$, were $\rec\fespace$ is the elliptic
reconstruction operator associated with the space $\fespace$, defined
in \ref{def:elliptic.reconstruction}.  In analogy with the fully
discrete notation in \S\ref{def:discrete.time}, we define the
\emph{semidiscrete elliptic reconstruction error} $\epsi:=\w-U$ and
the \emph{semidiscrete parabolic error} $\rho:=U-u$, keeping in mind
the warning \S\ref{sec:notational-warning}

We observe that while in the simplified semidiscrete setting one
assumes the discrete space $\fespace$ to be invariant in time, in the
fully discrete setting (cf. \S\ref{sec:duality.fully}) we will take
into account the possibility of the discrete space to change, with
respect to the timestep.  For instance, in an adaptive mesh refinement
scheme the space change derives from the mesh's modification from a
time to the next.

Correspondingly to the fully discrete case (\ref{eqn:error.equation}),
we may write the following \emph{semidiscrete the parabolic--elliptic error
  relation}:
\begin{equation}
    \label{eqn:error.equation.semidiscrete}
    \ltwop{\partial_t\rho(t)}\phi +\abil{\rho(t)}\phi
    =
    \ltwop{\partial_t\epsi(t)}\phi
    +\ltwop{\lproj\fespace f(t)-f(t)}\phi,
    \Foreach\phi\in\honezw,t\in(0,T].
\end{equation}
\subsection{The dual solution}
The concept of \emph{parabolic dual solution}, introduced first by
Eriksson \& Johnson~\cite{ErikssonJohnson:91} in the context of
\aposteriori error estimation, will be used now to obtain error
estimates out of (\ref{eqn:error.equation.semidiscrete}).  

For each $s\leq T$, consider the {\em dual solution} to be the
function
\begin{equation}
  z(x,t;s)=\zetas(x,t),\text{ for } x\in\W\text{ and }0\leq t\leq s,
\end{equation}
which satisfies $\zetas\in\leb2(0,T;\honezw)$,
$\partial_t\zetas\in\leb2(0,T;\sobh{-1}(\W))$, 
and solves the following backward parabolic {\em dual problem}:
\begin{equation}
  \label{eqn:dual.problem}
  \begin{split}
    -\ltwop{\partial_t\zetas(t)}\phi+\abil\phi{\zetas(t)}
    &=0,
    \Foreach\phi\in\honezw,t\in[0,s),
      \\
      \zetas(x,s)
      &=\rho(x,s)\Foreach x\in\W
  \end{split}
\end{equation}
for each $s\in[0,T]$.  Notice that $\phi$ can be taken to be time
dependent, with the appropriate differentiability properties.

The dual solution enjoys stability properties which we will use in the
sequel.  An immediate property is the usual energy identity
\begin{equation}
  \label{eqn:weak.stability}
  \Norm{\zetas(t)}^2+2\int_t^s\norm{\zetas}_a^2
  =\Norm{\rho(s)}^2\Foreach t\in[0,s].
\end{equation}
A more intricate stability property of $\zetas$ is given by the following result.
\begin{Lem}[Strong stability estimate {\cite[Lem. 4.2]{ErikssonJohnson:91}}]
\label{lem:strong.stability}

For each $s\in[0,T]$,
\begin{equation}
  \label{eqn:strong.stability}
  \qc{
  \int_0^s \Norm{\partial_t\zetas(t)}^2(s-t) \d t,
  \int_0^s \Norm{\contopera\zetas(t)}^2(s-t) \d t
  }
  \leq
  \frac14\Norm{\rho(s)}^2.
\end{equation}
\end{Lem}
\begin{Proof}
  For a fixed $s\in\clinter0T$, the change of variables
  \begin{equation}
    w(\vec x,t)=\zetas(\vec x,s-t)
  \end{equation}
  in the PDE (\ref{eqn:dual.problem}) implies
  \begin{equation}
    \ltwop{\pdt w(t)}\phi+\abil\phi{w(t)}=0
    \Foreach\phi\in\honezw,t\in\opclinter0T.
  \end{equation}
  Hence, testing with $\pdt w(t)t$ and integrating in time, we get
  \begin{equation}
    \begin{split}
      \int_0^s\Norm{\pdt w(t)}^2t\d t
      =&
      -\int_0^s\abil{\pdt w(t)}{w(t)}t\d t
      \\
      =&
      \int_0^s
      \frac12\norm{w(t)}_a^2
      -\frac12\pdt{\qb{\norm{w(t)}_a^2 t}}
      \d t
      \\
      \leq&
      \frac14\Norm{w(0)}^2
      -\frac12\qp{\Norm{w(s)}^2+\norm{w(s)}_a^2s},
    \end{split}
  \end{equation}
  where our last step relies standard energy identity:
  \begin{equation}
    \int_0^s\norm{w(t)}_a^2\d t
    =
    \frac12\Norm{w(0)}^2-\frac12\Norm{w(t)}^2.
  \end{equation}
  Thus we have
  \begin{equation}
    \int_0^s\Norm{\pdt w(t)}^2t\d t
    \leq
    \frac14\Norm{w(0)}^2,
  \end{equation}
  which is the first inequality in (\ref{eqn:strong.stability}); to
  obtain the second inequality, simply use the fact that $\pdt
  w(t)=\contopera w$.
\end{Proof}
\subsection{\Aposteriori error analysis via parabolic duality}
Integrating in \eqref{eqn:dual.problem} by parts in time implies that
\begin{equation}
  \begin{split}
    \ltwop{\rho(s)}{\phi(s)}
    &=\ltwop{\zetas(s)}{\phi(s)}
    \\
    &=\ltwop{\zetas(0)}{\phi(0)}
    +\int_0^s\ltwop{\partial_t\phi(t)}{\zetas(t)}
    +\abil{\phi(t)}{\zetas(t)}\d t,
  \end{split}
\end{equation}
for all $\phi\in\leb2(0,T;\honezw)$ such that
$\partial_t\phi\in\leb2(0,T;\sobh{-1}(\W))$.

Take $\phi=\rho$, use \eqref{eqn:error.equation.semidiscrete} and
assume $\projf{}-f=0$ momentarily---in the proof of Theorem
\ref{the:general.aposteriori.estimate} we shall remove this assumption---we
obtain
\begin{equation}
  \label{eqn:semidiscrete.basic.estimate}
  \Norm{\rho(s)}^2=\ltwop{\rho(0)}{\zetas(0)}
  +\int_0^s\ltwop{\partial_t\epsi(t)}{\zetas(t)}\d t.
\end{equation}
The first term on the right-hand side, is easily estimated, with
Lemma \ref{lem:strong.stability} in mind, as follows
\begin{equation}
  \ltwop{\rho(0)}{\zetas(0)}
  \leq\Norm{\rho(0)}\sup_{[0,s]}\Norm{\zetas}.
\end{equation}
As for the second term on the right-hand side of
\eqref{eqn:semidiscrete.basic.estimate} we have the choice of two
different ways for estimating it.
\begin{LetterList}
  \item
    \label{itm:direct.approach} 
    A direct estimate yields
    \begin{equation}
      \int_0^s\ltwop{\partial_t\epsi}{\zetas}
      \leq
      \sup_{[0,s]}\Norm{\zetas}\int_0^s\Norm{\partial_t\epsi}.
    \end{equation}
    Notice that the term $\partial_t\epsi$ can be estimated via
    elliptic \aposteriori error estimates because it is the difference
    between $\partial_t U$ and its reconstruction $\rec \partial_t
    U=\partial_t\rec U$.  Nonetheless a term involving $\partial_t\epsi$ is less
    desirable than one involving only $\epsi$.
  \item
    \label{itm:indirect.approach} 
    A less direct estimate, that would avoid the appearance of
    time derivatives in the indicator, is obtained by integrating by
    parts in time first
    \begin{equation}
      \int_0^s\ltwop{\partial_t\epsi}{\zetas}
      =
      \ltwop{\epsi(s)}{\zetas(s)}-\ltwop{\epsi(0)}{\zetas(0)}
      -\int_0^s\ltwop{\epsi(t)}{\partial_t\zetas(t)}\d t.
    \end{equation}
    The last integral can be then bounded as follows
    \begin{equation}
      \begin{split}
	\int_0^s{\epsi(t)}{\partial_t\zetas(t)}\d t
	&\leq
	\int_0^s\frac{\Norm{\epsi(t)}}{\sqrt{s-t}}
	\Norm{\partial_t\zetas(t)}\sqrt{s-t}\d t\\
	&\leq
	\left(\int_0^s\frac{\Norm{\epsi(t)}^2}{s-t}\d t\right)^{1/2}
	\left(\int_0^s\Norm{\partial_t\zetas(t)}^2(s-t)\d t\right)^{1/2}.
      \end{split}
    \end{equation}

    Unfortunately this bound turns out not to be useful, as it stands,
    due to the weight in the first integral on the last right-hand
    side. Namely, for this term to be finite it is necessary that
    $\epsi(t)=\oh(1)$ at $t=s$.  This means that the error between the
    discrete solution and its reconstruction should at least vanish at
    $s$.  Heuristically this can be interpreted as the mesh having to
    become infinitely fine as time gets closer to $s$: an unrealistic
    option.
\end{LetterList}

To circumvent this difficulty, without totally sacrificing
$\Norm{\epsi}$ to $\Norm{\partial_t\epsi}$, we compromise between
approach \ref{itm:direct.approach} and \ref{itm:indirect.approach} by
following through from \eqref{eqn:semidiscrete.basic.estimate} as
follows:  fix $r\in(0,s)$ (think of it as a close point to $s$), split
the integral and integrate by parts in time
\begin{equation}
  \label{eqn:semidiscrete.compromise}
  \begin{split}
    \Norm{\rho(s)}^2
    =&
    \ltwop{\zetas(0)}{\rho(0)}
    +\left(\int_0^r+\int_r^s\right)\ltwop{\partial_t\epsi}{\zetas}\\
    =&\ltwop{\zetas(0)}{\rho(0)-\epsi(0)}+\ltwop{\zetas(r)}{\epsi(r)}
    -\int_0^r\ltwop{\epsi}{\partial_t\zetas}
    +\int_r^s\ltwop{\partial_t\epsi}{\zetas}
    \\
    \leq& \sup_{[0,s]}\Norm{\zetas}\left(\Norm{e(0)}+\Norm{\epsi(r)}
    +\int_r^s\Norm{\partial_t\epsi}\right)\\
    &+\left(\int_0^r\Norm{\partial_t\zetas(t)}^2(s-t)\d t\right)^{1/2}
    \left(\int_0^r\frac{\Norm{\epsi(t)}^2}{s-t}\d t\right)^{1/2}.
  \end{split}
\end{equation}
The stability estimates \eqref{eqn:weak.stability} and
\eqref{eqn:strong.stability} imply that 
\begin{equation}
  \label{eqn:semidiscrete.estimate}
  \Norm{\rho(s)}\leq
  \Norm{e(0)}+\Norm{\epsi(r)}
    +\int_r^s\Norm{\partial_t\epsi}
    +\frac12\left(\int_0^r\frac{\Norm{\epsi(t)}^2}{s-t}\d t\right)^{1/2}.
\end{equation}
This discussion's outcome can be summarized into the following result.
\begin{The}[Semi-discrete duality-reconstruction \aposteriori error estimate]
  \label{the:semi-discrete-duality-reconstruction-estimate}
  Suppose that $f(t)\in\fespace$, for $t\in[0,T]$, and that there
  exists an \aposteriori elliptic $\leb2(\W)$-error estimator
  functional $\ef\cdot\fespace{\leb2(\W)}$, as defined in
  \S\ref{def:elliptic-error-estimator-functional},
  then the error occurring in the semi-discrete scheme
  \eqref{eqn:semidiscrete.scheme} obeys the \aposteriori bound
  \begin{equation}
    \begin{split}
      \sup_{t\in[0,s]}\Norm{U(t)-u(t)}
      \leq
      &\Norm{U(0)-u(0)}
      +L(s,r)\sup_{[0,s]}\ef {U}{\fespace}{\leb 2(\W)}\\
      &+(s-r)\sup_{[r,s]}\ef {\partial_t U}{\fespace}{\leb 2(\W)}
    \end{split}
  \end{equation}
  for each $s,r$, $0\leq s<r\leq T$, and with
  \begin{equation}
    \label{eqn:log.factor}
    L(s,r):=2+\frac12\sqrt{\log\frac s{s-r}}.
  \end{equation}
\end{The}
\begin{Proof}
  Fix $r$ and $s$ and use \eqref{eqn:semidiscrete.estimate} to get 
  \begin{equation}
    \label{eqn:proof:semidiscrete.estimate:start}
    \Norm{e(s)}
    \leq
    \Norm{e(0)}
    +\Norm{\epsi(r)}+\Norm{e(r)}
    +\int_r^s\Norm{\partial_t\epsi}
    +\frac12\psqrt{\int_0^r\frac{\Norm{\epsi(t)}^2}{s-t}\d t}.
  \end{equation}
  Basic manipulations and the use of the estimator functional
  $\ef\cdot\fespace{\leb2(\W)}$ leads to
  \begin{equation}
    \label{eqn:proof:semidiscrete.estimate:0-order-part}
    \begin{split}
      \Norm{\epsi(r)}
      &
      +\Norm{e(r)}
      +\frac12\psqrt{\int_0^r\frac{\Norm{\epsi(t)}^2}{s-t}\d t}
      \\
      &
      \leq
      \qp{2+\frac12\psqrt{\int_0^r\frac{\d t}{s-t}}}\sup_{0\leq t\leq s}\Norm{\epsi(t)}
      \\
      &
      \leq
      \qp{2+\frac12\sqrt{\log\frac s{s-r}}}
      \sup_{0\leq t\leq s}\ef{U(t)}\fespace{\leb2(\W)}
    \end{split}
  \end{equation}
  and
  \begin{equation}
    \label{eqn:proof:semidiscrete.estimate:1-order-part}
    \int_r^s\Norm{\partial_t\epsi}
    \leq
    (s-r)\sup_{r\leq t\leq s}\ef{\pdt U(t)}\fespace{\leb2(\W)}.
  \end{equation}
  The result follows by using
  (\ref{eqn:proof:semidiscrete.estimate:0-order-part}) and
  (\ref{eqn:proof:semidiscrete.estimate:1-order-part}) in
  (\ref{eqn:proof:semidiscrete.estimate:start}).
\end{Proof}
\begin{Cor}[Semi-discrete duality-residual \aposteriori estimates]
  \label{cor:semi-discrete-duality-residual-estimate}
If $\W$ is a convex domain in $\reals^d$ and $f(t)\in\leb2(\W)$ for
each $t\in[0,T]$, then the following \aposteriori error estimate holds
\begin{equation}
  \begin{split}
    \sup_{[0,s]}\Norm{U-u}\leq&\Norm{U(0)-u(0)}\\
    &+L(s,r)\sup_{[0,r]}\big(
    C_3\Norm{h^2({\pwA} - \disca{\fespace})U}
    \\
    &\phantom{+L(s,r)\sup_{[0,r]}\big(}
    +C_5\Norm{h^{3/2}J[U]}_\Sigma
    +C_7\Norm{h^2(\projf{}-f)}
    \big)\\
    &+(s-r)\sup_{[r,s]}\big(
    C_3\Norm{h^2({\pwA} - \disca{\fespace}) \partial_t U}
    \\
    &\phantom{+(s-r)\sup_{[r,s]}\big(}
    +C_5\Norm{h^{3/2}J[\partial_t U]}_\Sigma
    +\frac1{2\sqrt\alpha}\Norm{h(\projf{}-f)}
    \big).
  \end{split}
\end{equation}
\end{Cor}
\begin{Proof}
  From Theorem \ref{the:semi-discrete-duality-reconstruction-estimate}
  and Lemma \ref{lem:elliptic.aposteriori.estimates}, the result
  follows when $f(t)\in\fespace$ for all $t\in\opinter0T$.  To remove
  this assumption, 
\end{Proof}
\section{Estimates for the fully discrete scheme}
\label{sec:duality.fully}
Bearing in mind the techniques of the last section, we now turn our
attention to the analysis of the fully discrete scheme
\eqref{eqn:fully.discrete.scheme}.  For convenience, we switch
notation slightly and use the symbol $U$ (even without the superscript
$n$ sometime) for the fully discrete solution and its piecewise linear
interpolation now.  We introduce first some extra ``discrete-time''
notation to be used in this section.

\begin{Def}[duality time-accumulation coefficients]
\label{def:duality-time-coefficients}
In developing the error bounds via duality, we shall need the
following (logarithmic) {\em time accumulation coefficients}:
\begin{equation}
  \label{eqn:duality-time-accumulation-coefficients}
  \begin{gathered}
    b_n:=
    \begin{cases}
      \frac14\log\left(\frac{T-\tno}{T-\tn}\right),
      &
      \text{ for }\rangefromto n1{N-1},
      \\
      \frac18,
      &
      \text{ for }n=N,
    \end{cases}
    \\
    a_n
    :=
    \int_0^\TNO\frac{ l_n(t)\d t}{T-t}
    =
    \begin{cases}
      1-\lambda\qp{-\frac{\tau_1}T},
      &\text{ for }n=0,
      \\
      \lambda\left(\frac{\taun}{T-\tn}\right)
      -\lambda\left(-\frac{\tau_{n+1}}{T-\tn}\right),
      &\text{ for }\rangefromto n1{N-2},
      \\
      \lambda\left(\frac{\TAUNO}{\TAUN}\right)-1, 
      &\text{ for }n = N-1,
    \end{cases}
  \end{gathered}
\end{equation}
where
\begin{equation}
  \lambda(x):= 
  \begin{cases}
    (1+1/x)\log(1+x)
    &
    \text{ for }\norm{x}\in(0,1),
    \\
    1
    &
    \text{ for }x=0,
  \end{cases}
\end{equation}
which is an increasing function of $x$. We observe that the functions
$\lambda(x)-1$, $1-\lambda(-x)$ and $\lambda(x)-\lambda(-y)$ are
positive for $(x,y)\in(0,1)^2$, a fact that makes the coefficients
$a_n$ to be positive.  These coefficients can be appreciated
graphically in Figure~\ref{fig:time-accumulation-coefficients}.
\end{Def}
\begin{figure}[h]
  \caption{
    \label{fig:time-accumulation-coefficients}
    An example of the time accumulation coefficients $\seqi an$,
    $\seqi b n$ and $\seqi d n$ defined in
    (\ref{eqn:duality-time-accumulation-coefficients}) and
    (\ref{eqn:energy-time-accumulation-coefficients}),
    respectively. This is the situation for a uniform timestep $1/4$
    over the interval $\clinter0{10}$. All coefficients exhibit a
    backward decaying ``tail'' (cf. Theorem
    \ref{the:general.aposteriori.estimate}).  Noting how this tail is
    much heavier for $\seqi an$ and $\seqi b n$ than for $\seqi d n$ it
    follows that the energy estimator ``forgets'' much faster than the
    duality one.
 }  
  \newcommand{\xscale}{1}
  \newcommand{\yscale}{5}
  \begin{tikzpicture}[xscale=\xscale,yscale=\yscale]
  \foreach \xtick in {0,1,2,3,4,5,6,7,8,9,10}{
    \newcommand{\pointx}{\xtick}
    \draw (\xtick cm,0.025cm)--(\xtick cm,-0.025cm) 
    node[below]{$\pointx$};
  }
  \foreach \ytick in {0.0,0.25,0.5,0.75,1.0}{
    \newcommand{\pointy}{\ytick}
    \draw (0.1cm,\ytick cm)--(-0.1cm,\ytick cm)
    node[left]{$\pointy$};
  }
  
  \draw[-stealth] (0,0) -- (10.5,0) node[below]{$x$};
  \draw[-stealth] (0,0) -- (0,1.1) node[left]{$y$};
  
  
  \draw[color=c,line width=2pt] 
  plot file{Picture/an.table} node [c, below] {$a_n$};
  \draw[color=a,line width=2pt]
  plot file{Picture/bn.table} node [a, below] {$b_n$};
  \draw[color=b,line width=2pt]
  plot file{Picture/dn.table} node [b, above] {$d_n$};
\end{tikzpicture}    
\end{figure}
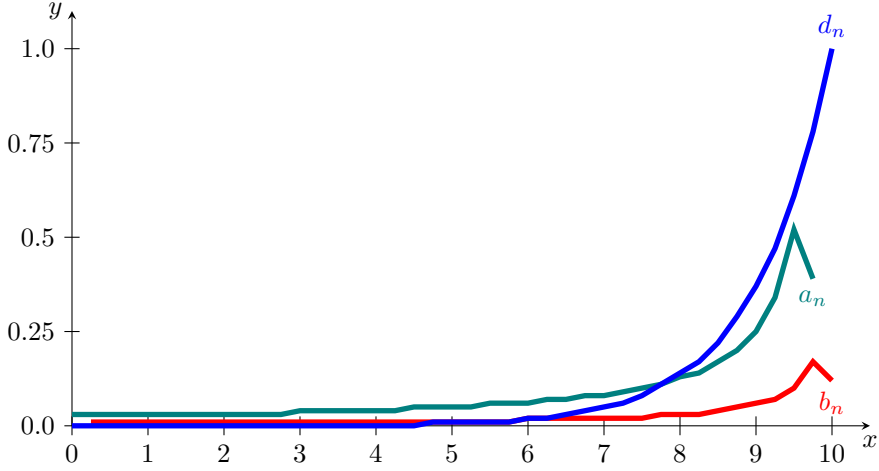
\begin{Lem}[duality time-accumulation coefficients properties]
  The coefficients $a_n$ and $b_n$, defined in
  \S\ref{def:duality-time-coefficients} for $\rangefromto n0N$,
  satisfy the following
  \begin{equation}
    \label{eqn:properties:duality-time-coefficients}
    \sumifromto n0{N-1}a_n=\log\frac T\TAUN
    ,
    \quad
    \int_{\tno}^\tn\frac{l_{n-1}(t)^2}{T-t}\d t
    \leq b_n
    \AND
    \sumifromto n1N b_n=\frac14\qp{\frac12+\log\frac T\TAUN}
  \end{equation}
\end{Lem}
\begin{Proof}
  The results follow from the definitions and basic calculus.
\end{Proof}
\begin{Def}[error indicators]
  \label{def:estimators}
  We suppose an \aposteriori elliptic error estimator functional
  $\ef\cdot\cdot\cdot$, as defined in
  \S\ref{def:elliptic-error-estimator-functional}, is available and we
  introduce the following (time-local) {\em $\cE$-based spatial error
    indicators}
  \begin{gather}
    \label{eqn:def:fully-discrete-reconstruction-indicator}
    \indrec n:=\ef{U^n}{\fez n}{\leb2(\W)},
    \\
    \label{eqn:def:fully-discrete-spatial-parabolic-indicator}
    \indspace_n:={\taun}\ef{\bdisc\un}{\fez n\cap\fez {n-1}}{\leb2(\W)},
  \end{gather}
  and the {\em time error indicator}
  \begin{equation}
    \label{eqn:def:time-indicator}
    \indtime n:=\Norm{\ano\uno-\an\un}=
    \begin{cases}
      \frac12\Norm{\projf 1-\discn U^1-A^0U^0}
      &
      \text{ for }n=1,
      \\
      \frac12\Norm{\bdisc\left(\projf n-\discn\un\right)}\taun
      &
      \text{ for }\rangefromto n2N.
    \end{cases}
  \end{equation}
  or, in some cases, the alternative version
  \begin{equation}
    \label{eqn:def:alternative-time-indicator}
    \indtime n:= \Norm{\uno-\un}+\indspace_n.
  \end{equation}
  In the numerical experiments we only use definition (\ref{eqn:def:time-indicator})
  $\indtime n$.
    
  We also introduce the {\em data approximation error indicator}
  \begin{equation}
    \label{eqn:def:data-approximation-indicator}
    \inddata n:=\int_\tno^\tn\Norm{f^n-f(t)}\d t,
  \end{equation}
  the associated \emph{global data approximation indicator}
  \begin{equation}
    \tilde\beta_N:=
    \begin{cases}
      \sum_{n=1}^N\inddata n
      &\text{ if }
      \f n=f(\tn)
      ,
      \\
      \inddata N
      +2
      \psqrt{
        \sum_{n=1}^{N-1}b_n\inddata n^2
      }
      &\text{ if }
      \f{n}=\int_\tno^\tn{f(t)\d{t}}/\taun,
    \end{cases}
  \end{equation}
  and the \emph{mesh change error indicator function}
  \begin{equation}
    \label{eqn:def:coarsening-indicator}
    \indcoarse n := \frac{\lproj n\uno-\uno}\taun+\projf n-\fn
    =(\lproj n-\Id)(\taunm\uno+\fn).
  \end{equation}
\end{Def}
\begin{Obs}[smooth data approximation]
  For $f$ smooth enough, we can redefine the indicator $\inddata n$ in
  relation (\ref{eqn:def:data-approximation-indicator}) by the right
  hand side of the following inequality
  \begin{equation}
    \int_\tno^\tn\Norm{f^n-f(t)}\d t
    \leq\Norm{\partial_t f}_{\leb1(\In,\leb2(\W))}\taun.
  \end{equation}
\end{Obs}
\begin{The}[general duality \aposteriori parabolic-error estimate]
  \label{the:general.aposteriori.estimate}
  Let $u$ be the exact solution of (\ref{eqn:continuous.heat}),
  $\seqs{U^n}{\rangefromto n0N}$ the (corresponding fully discrete)
  solution of (\ref{eqn:fully.discrete.scheme}) and $\w^n=\rec n\un$
  the elliptic reconstruction of $\un$, for $\rangefromto n0N$ as
  defined by (\ref{eqn:reconstruction}). 
  Then, with reference to Definition \ref{def:estimators}, the following
  \aposteriori error estimate holds
  \begin{equation}
    \label{eqn:general-aposteriori-estimate}
    \begin{split}
      \Norm{\w^N-u(T)}
      \leq
      &\Norm{U^0-u(0)}
      +
      \psqrt{
        \sum_{n=0}^{N-1}a_n\indrec n^2
      }
      +\indspace_N
      +\psqrt{\sum_{n=1}^Nb_n\indtime n^2}
      \\&
      +\sqrt{\frac{\TAUN}2}\Norm{\indcoarse N h_N}
      +\psqrt{\sum_{n=1}^{N-1}b_n\Norm{\indcoarse n\hn^2}^2}
      +\tilde\beta_N.
    \end{split}
  \end{equation}
\end{The}
The proof of this result is the object of \S\ref{sec:proof-duality}.
We state and application of this result, which we will also prove
later in \S~\ref{sec:proof-corollary-fully-discrete}.
\begin{Cor}[Duality \aposteriori full error estimates]
  \label{cor:fully-discrete-estimate}
  With the same notation as in \ref{the:general.aposteriori.estimate}
  and supposing that $\f{n}=\int_\tno^\tn{f(t)\d{t}}/\taun$ we have 
  \newcommand{\dummy}{+\sqrt{1+\log\frac T{\TAUN}}\bigg(}
      \begin{multline}
        \Norm{\UN-u(T)}
        \leq
        \Norm{\UZERO-u(0)}
        +\sqrt{\frac\TAUN2}\Norm{\indcoarse N h_N}
        +\sum_{n=1}^N\tau_n\Norm{\partial_t f}_{\leb1(I_n;\leb2(\W))}  
        \\
        \dummy
        \max_{\rangefromto n0N}\indrec n
        +2
        \max_{\rangefromto n1{N-1}}\Norm{\indcoarse n h_n^2}
        +\frac12
        \max_{\rangefromto n1N}
        \indtime n   
        \bigg)
        .
      \end{multline}
\end{Cor}
\begin{Obs}[comparison between Theorem 
    \ref{the:general.aposteriori.estimate} and Corollary
    \ref{cor:fully-discrete-estimate}] Corollary
  \ref{cor:fully-discrete-estimate} has a simpler estimate than
  Theorem \ref{the:general.aposteriori.estimate} in that it involves
  less terms and does not require as much memory.  Notice however,
  that from an error bound view-point, the Theorem's tighter bound may
  be more effective as the time accumulation is not as strict as in
  the Corollary.  This is especially true in problems, typical in the
  parabolic setting, where the initial error may be very big and gets
  damped with time.
\end{Obs}
\section{Proof of Theorem \ref{the:general.aposteriori.estimate}}
\label{sec:proof-main-results}
\label{sec:proof-duality}
As with the semidiscrete case that we dealt with in
\secref{sec:duality.semi} to prove
Theorem~\ref{the:semi-discrete-duality-reconstruction-estimate}, our
starting point to prove (\ref{eqn:general-aposteriori-estimate}) is
the fully discrete analog of \eqref{eqn:semidiscrete.basic.estimate},
which is readily obtained from \eqref{eqn:error.equation} and
\eqref{eqn:dual.problem}:
\begin{multline}
  \label{eqn:discrete.basic.estimate}
  \Norm{\rho(T)}^2
  =
  \ltwop{\rho(0)}{\zetat(0)}
  +\int_0^T\ltwop{\partial_t\epsi(t)}{\zetat(t)}\d t\\
  \!+\!\!\sum_{n=1}^N\!
  \int_\tno^\tn\!\!\!\!\!\abil{\w(t)-\wn}{\zetat(t)\!}
  \!
  +\ltwop{\indcoarse n}{\zetat(t)}+\ltwop{\fn-f(t)}{\zetat(t)}
  \d t.
\end{multline}
We recall that $\rho$ and $\epsi$ are defined in  the functions $\indcoarse n$ are defined, for each $\rangefromto n1N$,
by (\ref{eqn:def:coarsening-indicator}) and $\zetat$ is the solution of the
dual problem~(\ref{eqn:dual.problem}) with $s=T$.
\subsection{Space error estimate}
The first two terms are estimated, similarly to
\eqref{eqn:semidiscrete.compromise}, as follows
\begin{multline}
  \ltwop{\rho(0)}{\zetat(0)}+\int_0^T\ltwop{\partial_t\epsi(t)}{\zetat(t)}\d
  t
  \leq
  \\ 
  \Norm{\rho(T)}
  \qp{
    \Norm{e(0)}+\Norm{\epsi(\TNO)}
    +\int_\TNO^T\Norm{\partial_t\epsi}
    +\frac12
    \psqrt{\int_0^\TNO\frac{\Norm{\epsi(t)}^2}{T-t}\d t}
  }.
\end{multline}
To proceed we observe that
\begin{equation}
  \label{eqn:proof:space-error-estimate:higher-order-terms}
  \int_\TNO^T\Norm{\partial_t\epsi} 
  =\Norm{\epsi^N-\epsi^{N-1}}\leq\indspace_N,
\end{equation}
and, by convexity and affinity of $l_n$, this implies
\begin{equation}
  \label{eqn:estimate.reconstruction.part}
  \begin{split}
    \int_0^\TNO\frac{\Norm{\epsi(t)}^2}{T-t}\d t
    =
    &\int_0^\TNO\frac{\Norm{\sum_{n=0}^{N-1}
	\epsi^n l_n(t)}^2}{T-t}\d t
    \\ 
    \leq
    &\sum_{n=0}^{N-1}\Norm{\epsi^n}^2
    \int_0^\TNO\frac{l_n(t)}{T-t}\d t
    =
    \sum_{n=0}^{N-1}a_n\indrec n^2.
    \end{split}
\end{equation}

Thus we obtain
\begin{multline}
  \label{eqn:estimate.space}
  \ltwop{\rho(0)}{\zetat(0)}+\int_0^T\ltwop{\partial_t\epsi(t)}{\zetat(t)}\d
  t
  \\
  \leq
  \Norm{\rho(T)}
  \qpBigg{
  \Norm{e(0)}
  +
  \psqrt{
    \sum_{n=0}^{N-1}a_n\indrec n^2
  }
  +\indspace_N
  }.
\end{multline}
\subsection{Time error estimate}
The third term in (\ref{eqn:discrete.basic.estimate}), which accounts
mainly for the time error, can be bounded as follows
\begin{equation}
  \begin{split}
    &\sum_{n=1}^N 
    \int_\tno^\tn
    \abil{\w(t)-\wn}{\zetat(t)}\d t
    \\
    &\leq 
    \sum_{n=1}^N\int_\tno^\tn\Norm{\wno-\wn}l_{n-1}(t)
    \Norm{\contopera\zetat(t)}\d t
    \\
    &
    \leq \frac12\Norm{\rho(T)}
    \psqrt{
      \sum_{n=1}^N
      \Norm{\wno-\wn}^2
      \int_\tno^\tn\frac{l_{n-1}(t)^2}{T-t}\d t
    }
    \\
    &
    \leq 
    \frac12\Norm{\rho(T)}
    \psqrt{
      \frac12\Norm{\w^N-\w^{N-1}}^2
      +\sum_{n=1}^{N-1}\Norm{\wn-\wno}^2
      \log\left(\frac {T-\tno}{T-\tn}\right)
    }
    .
  \end{split}
\end{equation}
Notice that how in the time integral on the last interval $(\TNO,\TN]$
    the numerator $l_{N-1}(t)=\Oh(T-t)$ compensates for the
    singularity of $1/(T-t)$.  The terms $\Norm{\wno-\wn}$ appearing
    in this estimate still need to be estimated, as there is no
    explicit knowledge of the reconstructed functions $\wn=\rec n
    U^n$.  These terms can be dealt with in two different ways.
\begin{LetterList}
\item
  One way to estimate these terms is given by:
  \begin{equation}
    \begin{split}
      \Norm{\wno-\wn}
      &\leq\Norm{\uno-\un}+\Norm{\wno-\wn-\uno+\un}\\
      &=\Norm{\uno-\un}+\taun\Norm{\partial_t\epsi^n}
      \leq\Norm{\uno-\un}+\indspace_n,
    \end{split}
  \end{equation}
  for all $t\in I_n$.
  Thus we obtain the estimate
  \begin{equation}
    \label{eqn:estimate.time.a}
    \begin{split}
      \sum_{n=1}^N 
      \int_\tno^\tn
      \!\!\!\!
      \abil{\w(t)-\wn}{\zetat(t)}\d t
      &
      \leq
      \Norm{\rho(T)}
      \psqrt{
        \sum_{n=1}^Nb_n
        \powqp2{
          \Norm{\uno-\un}
          +\indspace_n
        }
      }
      \\
      &
      =
      \Norm{\rho(T)}
      \psqrt{
        \sum_{n=1}^Nb_n\indtime n^2
      }
      ,
    \end{split}
  \end{equation}
  by using $\indtime n$'s alternative definition
  (\ref{eqn:def:alternative-time-indicator}).
\item 
  Another way to estimate $\Norm{\wno-\wn}$ consists in using again the
  definition of elliptic reconstruction and the PoincarÃ© inequality as
  follows:
  \begin{equation}
    \begin{split}
      \Norm{\wno-\wn}^2
      &\leq \constext{PF}\abil{\wno-\wn}{\wno-\wn}
      \\
      &= \constext{PF}\ltwop{\ano\uno-\an\un}{\wno-\wn}
      \\
      &\leq \constext{PF}\Norm{\ano\uno-\an\un}\Norm{\wno-\wn},
    \end{split}
  \end{equation}
  thus obtaining
  \begin{equation}
    \Norm{\wno-\wn}\leq \constext{PF}\Norm{\ano\uno-\an\un}
  \end{equation}
  Using the definition of $\indtime n$ in (\ref{eqn:def:time-indicator}) yields
  \begin{multline}
    \label{eqn:time.estimate.b}
    \sum_{n=1}^N 
    \int_\tno^\tn\abil{\w(t)-\wn}{\zetat(t)}\d t
    \\
    \leq
    \Norm{\rho(T)}
    \psqrt{
      \sum_{n=1}^Nb_n 
      \constext{PF}^2
      \qp{\ano\uno-\an\un}^2}
    =
    \Norm{\rho(T)}
    \psqrt{
      \sum_{n=1}^Nb_n\theta_n^2
    }
    .
  \end{multline}
\end{LetterList}
\subsection{Mesh-change error estimates}
To bound the third term in \eqref{eqn:discrete.basic.estimate}, we use
$\indcoarse n$'s $\leb2(\W)$ orthogonality and the orthogonal
projector $\funk{\Pi^n}{\leb2(\W)}{\fes n}$ as follows:
\begin{equation}
  \begin{split}
      \sum_{n=1}^N&\int_\tno^\tn
      \ltwop{\indcoarse n}{\zetat(t)}\d t
      =
      \sum_{n=1}^{N}\int_\tno^\tn
      \ltwop{\indcoarse n}
	    {\zetat(t)-\Pi^n\zetat(t)}\d t
      \\
      &\leq
      \sum_{n=1}^{N-1}\int_\tno^\tn
      \Norm{\indcoarse n\hn^2}
      \norm{\zetat(t)}_2\d t
      +\int_\TNO^\TN\Norm{\indcoarse N h_N}
      \norm{\zetat(t)}_1\d t
      \\
      &\leq\frac12\Norm{\rho(T)}\Bigg(
      \sqrt{2\TAUN}\Norm{\indcoarse N h_N}
      +\left(\sum_{n=1}^{N-1}
      \Norm{\indcoarse n\hn^2}^2
      \log\left(\frac{T-\tno}{T-\tn}\right)
      \right)^{1/2}\Bigg)
      \\
      &=\Norm{\rho(T)}
      \left(\sqrt{\frac{\TAUN}2}\Norm{\indcoarse N h_N}+
      \left(\sum_{n=1}^{N-1}b_n\Norm{\indcoarse n\hn^2}^2\right)^{1/2}
      \right).
  \end{split}
  \label{eqn:estimate.data}
\end{equation}
Here we have used the fact that $\W$ is convex in order to apply the
estimate
\begin{equation}
  \norm{\zetat(t)}_2\leq\Norm{\contopera\zetat(t)},
\end{equation}
and then apply the strong stability estimate
\eqref{eqn:strong.stability}.
\subsection{Data-approximation error estimates}
The fourth term in \eqref{eqn:discrete.basic.estimate} can be bounded
in two different ways depending on which definition for $\f n$
appearing in the fully discrete scheme
(\ref{eqn:fully.discrete.scheme}) is chosen.

\begin{LetterList}
\item If $\fn=f^n$ then we can proceed as follows
  \begin{equation}
    \label{eqn:estimate.data.time}
    \begin{split}
      \sum_{n=1}^N\int_\tno^\tn\ltwop{\fn-f(t)}{\zetat(t)}\d t
      &\leq
      \sum_{n=1}^N\max_{I_n}\Norm{\zetat}\int_\tno^\tn\Norm{\fn-f(t)}\d
      t\\
      &\leq
      \Norm{\rho(T)}\sum_{n=1}^N\inddata n.
    \end{split}
  \end{equation}
\item If instead of $\f n=f(\tn)$ we have
  $\f{n}=\int_\tno^\tn{f(t)\d{t}}/\taun$, which is the $\leb2$
  projection of $f$ onto constants in time, then we can exploit the
  orthogonality and write, for each $\rangefromto n1{N-1}$
  \begin{equation}
    \begin{split}
      \int_\tno^\tn\ltwop{\fn-f(t)}{\zetat(t)}\d t
      =
      \int_\tno^\tn\ltwop{\fn-f(t)}{\zetat(t)-\zetat(t_{n-1})}\d t
      \\
      \leq\max_{t\in\In}\Norm{\zetat(t)-\zetat(\tno)}
      \int_\tno^\tn\Norm{\fn-f(t)}\d t.
    \end{split}
  \end{equation}
  By noticing that
  \begin{equation}
    \begin{split}
      \max_{t\in\In}\Norm{\zetat(t)-\zetat(\tno)}
      &=
      \max_{t\in\In}\Norm{\int_{\tno}^t\partial_s\zetat(s)\d s}
      \\
      &\leq\max_{t\in\In}\int_{\tno}^t\Norm{\partial_s\zetat(s)}\d s
      \leq\int_{\tno}^{\tn}\Norm{\partial_t\zetat}
      \\
      &\leq
      \log\left(\frac{T-\tno}{T-\tn}\right)^{1/2}
      \psqrt{
        \int_\tno^\tn\Norm{\partial_t\zetat(t)}^2(T-t)\d t
      }
      \\
      &=2b_n^{1/2}
      \psqrt{
        \int_\tno^\tn\Norm{\partial_t\zetat(t)}^2(T-t)
        \d t}
      .
     \end{split}
  \end{equation}
Summing up, applying \CBS, and
using the strong stability estimate \eqref{eqn:strong.stability} we
obtain
\begin{equation}
  \label{eqn:estimate.data.time-improved}
  \sum_{n=1}^N\int_\tno^\tn\ltwop{\fn-f(t)}{\zetat(t)}\d t
  \leq
  \Norm{\rho(T)}
  \qp{
    \inddata N
    +2
    \psqrt{
      \sum_{n=1}^{N-1}b_n\inddata n^2
    }
  }
  .
\end{equation}

Employing estimates \eqref{eqn:estimate.space},
(\ref{eqn:estimate.time.a}) (or (\ref{eqn:time.estimate.b})),
\eqref{eqn:estimate.data} and (\ref{eqn:estimate.data.time}) (or
\eqref{eqn:estimate.data.time-improved}) into the relation
\eqref{eqn:discrete.basic.estimate} we obtain the result of Theorem
\ref{the:general.aposteriori.estimate}.\hfill\qed
\subsection{Proof of Corollary~\ref{cor:fully-discrete-estimate}} 
\label{sec:proof-corollary-fully-discrete}
Referring to the notation introduced in \S~\ref{def:estimators}, the
indicator $\indspace_n$ defined by
(\ref{eqn:def:fully-discrete-spatial-parabolic-indicator}) and
appearing in (\ref{eqn:general-aposteriori-estimate}) can be
substituted by the more ``practical'' one: $\indrec n+\indrec{n-1}$.
To see this let us first revisit
estimate~(\ref{eqn:proof:space-error-estimate:higher-order-terms}) and
recall definition
(\ref{eqn:def:fully-discrete-reconstruction-indicator}) to write
\begin{equation}
  \Norm{\pdt\epsi^N}\taun=\Norm{\epsi^N-\epsi^{N-1}}
  \leq\Norm{\epsi^N}+\Norm{\epsi^{N-1}}\leq\indrec N+\indrec {N-1}.
\end{equation}
It follows, that
\begin{multline}
  \Norm{\w^N-u(T)}
  \leq
  \Norm{U^0-u(0)}
  +
  \psqrt{
    \sum_{n=0}^{N-1}a_n\indrec n^2
  }
  +\indrec{N-1}+\indrec N
  \\
  +\psqrt{\sum_{n=1}^Nb_n\indtime n^2}
  +\sqrt{\frac{\TAUN}2}\Norm{\indcoarse N h_N}
  +\psqrt{\sum_{n=1}^{N-1}b_n\Norm{\indcoarse n\hn^2}^2}
  +\tilde\beta_N
  .
\end{multline}
To close use the splitting $e=\rho-\epsi$ to obtain
\begin{equation}
  \begin{split}
    \Norm{e^N}
    \leq
    &\Norm{e^0}
    +
    \psqrt{
      \sum_{n=0}^{N-1}a_n\indrec n^2
    }
    +\indrec{N-1}+\indrec N
    \\&
    +\psqrt{\sum_{n=1}^Nb_n\indtime n^2}
    +\sqrt{\frac{\TAUN}2}\Norm{\indcoarse N h_N}
    \\&
    +\psqrt{\sum_{n=1}^{N-1}b_n\Norm{\indcoarse n\hn^2}^2}
    +  \inddata N
      +2
      \psqrt{
        \sum_{n=1}^{N-1}b_n\inddata n^2
      }
  \end{split}
\end{equation}

We notice that the former estimate implies the more traditional one
\cite{ErikssonJohnson:91}
\begin{equation}
  \label{eqn:estimate.space.max}
  \begin{split}
    \int_0^\TNO\frac{\Norm{\epsi(t)}^2}{T-t}\d t
    \leq\max_{\rangefromto n0{N-1}}\Norm{\epsi^n}^2
    \left(1+\sum_{n=1}^{N-1}a_n\right)\\
    =\left(1+4L(T,\TNO\right)^2)\max_{\rangefromto n0{N-1}}\Norm{\epsi^n}^2
  \end{split}
\end{equation}  
where $L(T,\TNO)$ is the logarithmic factor defined in
\eqref{eqn:log.factor}.

Also here, the indicator can be simplified if we relax the bound by
taking the maximum norm in time as follows:
\begin{multline}
  \label{eqn:estimate.time.max}
  \sum_{n=1}^N \int_\tno^\tn\abil{\w(t)-\wn}{\zetat(t)}\d t
  \\
  \leq
  \Norm{\rho(T)}\max_{\rangefromto n1N}\indtime n
  \sqrt{\left(\frac18+L(T,t_{N-1})^2\right)}.
\end{multline}

  As with the space and time estimates,
this estimate can be simplified, with some loss of sharpness, as
follows
\begin{multline}
  \label{eqn:estimate.data.max}
      \sum_{n=1}^N\int_\tno^\tn
      \ltwop{\indcoarse n}{\zetat(t)}\d t
      \\
      \leq
      \Norm{\rho(T)}
      \left(\sqrt{\frac\TAUN2}\Norm{\indcoarse N h_N}
      +L(\TN,\TNO)
      \max_{\rangefromto n1{N-1}}\Norm{\indcoarse n h_n^2}
      \right).
\end{multline}

Like earlier estimates, this estimate can be further simplified, by
taking the maximum and slightly relaxing it, into
\begin{multline}
  \sum_{n=1}^N\int_\tno^\tn\ltwop{\fn-f(t)}{\zetat(t)}\d t
  \\
  \leq
  \Norm{\rho(T)}\left(\inddata N
  +2L(\TN,\TNO)
  \max_{\rangefromto n1{N-1}}\inddata n\right).
\end{multline}
\end{LetterList}
\hfill\qed
\section{The energy-reconstructive approach}
\label{sec:energy}
In a previous paper~\cite{lakkis-makridakis:06}, we analyzed the
combination of classical energy methods for parabolic equations with
the elliptic reconstruction to obtain \aposteriori
$\leb\infty(\leb2(\W)$-error estimates.  In this section we give a
similar analysis that yields tighter bounds with respect to time
accumulation, that will be compared with the ones derived by duality.
\begin{The}[Semi-discrete energy-reconstruction \aposteriori error estimate]
  \label{the:semi-discrete-energy-reconstruction-estimate}
  Let the notation and conditions of Theorem
  \ref{the:semi-discrete-duality-reconstruction-estimate} hold then
  the following \aposteriori bound is true
  \begin{equation}
    \begin{split}
      \sup_{[0,s]}\Norm{U-u}
      &\leq
      \Norm{U(0) - u(0)} 
      +
      \ef{U(0)}{\fes{}}{\leb 2(\W)}
      +
      \sup_{[0, s]} \ef{U}{\fes{}}{\leb 2(\W)}
      \\
      &\phantom{  } +
      \int_0^s \exp{\left(\frac{\alpha}{c_0^2}(s-t)\right)}
      \ef{\pd t U}{\fes{}}{\leb 2(\W)} \d t.
    \end{split}
  \end{equation}
\end{The}
\begin{Proof}
  Testing equation (\ref{eqn:error.equation.semidiscrete}) with $\rho$
  and noting $f(t)\in\fes{}$ yields
  \begin{equation}
    \frac{1}{2}\dt\Norm{\rho}^2
    +
    \norm{\rho}_a^2 
    =
    \ltwop{\pd t \epsilon}{\rho}
  \end{equation}
  In view of the PoincarÃ©--Friedrichs inequality and the equivalence
  of the energy norm to the \sobh1{(\W)} semi norm
  (\ref{eqn:one.energy.inequality})
  \begin{equation}
    \frac{1}{2}\dt\Norm{\rho}^2
    +
    \norm{\rho}_a^2 
    \geq
    \frac{1}{2}\dt\Norm{\rho}^2
    +
    \frac{\alpha}{c_0^2}\Norm{\rho}^2,
  \end{equation}
  and hence
  \begin{equation}
    \begin{split}
      \frac{1}{2}\dt\Norm{\rho}^2
      +
      \frac{\alpha}{c_0^2}\Norm{\rho}^2
      &\leq
      \ltwop{\pd t \epsilon}{\rho}
      \\
      &\leq \Norm{\rho} 
      \Norm{\pd t \epsilon} 
      .
    \end{split}
  \end{equation}
  Dividing through by $\Norm{\rho}$ gives
  \begin{equation}
    \dt \Norm{\rho} 
    +
    \frac{\alpha}{c_0^2}\Norm{\rho}
    \leq
    \Norm{\pd t \epsilon} 
    .
  \end{equation}
  Using the integrating factor
  $\exp\qp{\fracl{\alpha t}{c_0^2}}$ we conclude that
  \begin{equation}
    \label{eq:semi-discrete-error-equation}
    \Norm{\rho(t)} 
    \leq 
    \Norm{\rho(0)}
    +
    \int_0^t
    \exp{\left(\frac{\alpha}{c_0^2}(s - t)\right)}
    \Norm{\pd t \epsilon(s)} 
    \d s.
  \end{equation}
\end{Proof}
\begin{Cor}[Semi-discrete energy-residual \aposteriori estimates]
  Let the assumptions of Corollary
  \ref{cor:semi-discrete-duality-residual-estimate} hold, then the
  following \aposteriori bound holds
  \begin{equation}
    \begin{split}
      \sup_{[0,s]}\Norm{U-u}
      &\leq
      \Norm{U(0)-u(0)} 
      +
      C_3\Norm{h^2(\pwA - \disca{\fespace})U(0)} 
      +
      C_5\Norm{h^{3/2}J[U(0)]}_\Sigma 
      \\
      &\phantom{\leq}+\sup_{[0,s]}
      \qp{ 
        C_3\Norm{h^2({\pwA} - \disca{\fespace})U} 
        +
        C_5\Norm{h^{3/2}J[U]}_\Sigma 
      }
      \\
      &\phantom{\leq}
      +
      \int_0^s \exp{\left(\frac{\alpha}{c_0^2}(s-t)\right)}
      \sup_{[0,s]}
      \qpbigg{
        C_3\Norm{h^2({\pwA} - \disca{\fespace})\partial_t U} 
        \\
        &\phantom{\leq \exp{\left(\frac{\alpha}{c_0^2}(s-t)\right)}}
        +
        C_5\Norm{h^{3/2}J[\partial_t U]}_\Sigma
        +
        \Norm{\projf{}-f}
      }
      \d t.
    \end{split}
  \end{equation}
\end{Cor}
\begin{Proof}
  Removing the assumption $f(t)\in\fes{}$ from the proof of Theorem
  \ref{the:semi-discrete-energy-reconstruction-estimate} gives
  \begin{equation}
    \Norm{\rho(t)} 
    \leq 
    \Norm{\rho(0)}
    +
    \int_0^t
    \exp{\left(\frac{\alpha}{c_0^2}(s - t)\right)}
    \left( 
    \Norm{\pd t \epsilon(s)} 
    + \Norm{(\projf{} - f)(s)}
    \right)
    \d s
  \end{equation}
  as an analog of (\ref{eq:semi-discrete-error-equation}). The
  splitting $e(t) = \rho(t) - \epsilon(t)$ together with the error
  estimates
  from Lemma \ref{lem:elliptic.aposteriori.estimates} yield the
  desired result.
\end{Proof}
\begin{Def}[energy time accumulation coefficients]
  \label{def:energy-accumulation-coefficients}
  The \emph{energy time-ac\-cu\-mu\-lation function} is defined as
  \begin{equation}
    \label{eqn:energy-time-accumulation-function}
    d(t,s):=
    \exp\qp{a(t-s)}
    ,
    \quad
    0\leq t<s\leq T
    ,
    \quad
    a:=\fracl\alpha{c_0^2}
  \end{equation}
  where $c_0$ is the PoincarÃ©--Friedrichs constant and $\alpha$ is the
  coercivity constant (\ref{eqn:bounded.bilinear}); we denote $d(t,T)=:d(t)$. The
  \emph{energy time-accumulation coefficients} are defined, for $0\leq
  n<m\leq N$, by
  \begin{equation}
    \label{eqn:energy-time-accumulation-coefficients}
    d^m_n
    :=\int_\tno^\tn d(t,t_m)\d t
    =\frac1a\exp\qp{a(\tn-t_m)}\qp{1-\exp(-a\taun)}.
  \end{equation}
  When $m=N$ we drop it and simply write $d_n$ instead of $d^m_n$.
  Note the useful recursive relation
  \begin{equation}
    \label{eqn:energy-time-accumulation-coefficients:recursion}
    d^{m+1}_n=d^m_n\exp(a\tau_{m+1}).
  \end{equation}
\end{Def}
\begin{The}[general energy \aposteriori parabolic-error estimate]
  \label{the:general.energy.aposteriori.estimate}
  Making use of the same notation as in Theorem
  \ref{the:general.aposteriori.estimate} the following \aposteriori
  estimate holds
  \begin{equation}
    \label{eqn:energy-aposteriori-estimate}
    \begin{split}
      \max_{t_n\in[0,T]} \Norm{U^n - u(t_n)} 
      &
      \leq
      \Norm{\rho(0)}
      +
      \max_{n\in[0:N]} \varepsilon_n
      +
      2\sum_{n=1}^N
      \qp{\eta_n + \beta_n + \gamma_n +\theta_n}d_n
      \\
      &
      =:
      \cE_0+\cE_\infty(N)+\cE_1(N).
    \end{split}
  \end{equation}
\end{The}
\begin{Obs}[timestepping the error estimate in practice]
  In practice, the bound (\ref{eqn:energy-aposteriori-estimate}) has
  to be used at each ``final'' time $t_m$ instead of $T$.  When
  stepping from time $t_{m-1}$, to the next one, say $t_m$ then,
  thanks to the recursion
  (\ref{eqn:energy-time-accumulation-coefficients:recursion}), it is
  straightforward to update the new error estimator:
  \begin{gather}
    \cE_\infty(m)=\maxofset{\cE(m-1),\ep_m}
    \intertext{and}
    \cE_1(m)=\cE_1(m-1)\exp\qp{a\tau_m}
    +d^m_{m-1}\qp{\eta_n + \beta_n + \gamma_n +\theta_n}.
  \end{gather}
  This is an advantage of using the energy estimate
  (\ref{eqn:energy-aposteriori-estimate}) as an alternative to
  (\ref{eqn:general-aposteriori-estimate}), where the indicators
  have to be stored for all time-steps and the sums recomputed at each
  timestep.
\end{Obs}
\begin{proof}
\label{sec:proof-energy}
We utilize the arguments of \cite{lakkis-makridakis:06} under the approach
described in Theorem
\ref{the:semi-discrete-energy-reconstruction-estimate}. The starting
point for this estimate is the parabolic error identity
(\ref{eqn:error.equation}) tested with $\rho$ as follows
\begin{equation}
  \begin{split}
    \frac{1}{2}\dt\Norm{\rho}^2 + \Norm{\rho}_a^2
    &=
    \ltwop{\partial_t\epsi}\rho
    +\abil{\w-\wn}\rho
    +\taunm\ltwop{\lproj n\uno-\uno}\rho\\
    &\phantom{= \ltwop{\partial_t\epsi}\rho
    +\abil{\w-\wn}\rho}
    +\ltwop{\projf n-f}\rho\\
    &=:
    \cI_1 + \cI_2
    + \cI_3 + \cI_4.
  \end{split}
\end{equation}
Analogously to the semidiscrete we make use of a PoincarÃ©--Friedrichs
inequality and the coercivity of $a$ to absorb the energy norm into
the $\leb{2}(\W)$ norm as follows
\begin{equation}
  \frac{1}{2}\dt\Norm{\rho}^2
    +
    \norm{\rho}_a^2 
    \geq
    \frac{1}{2}\dt\Norm{\rho}^2
    +
    \frac{\alpha}{c_0^2}\Norm{\rho}^2.
\end{equation}
Giving
\begin{equation}
  \begin{split}
    \frac{1}{2}\dt\Norm{\rho}^2
    +
    \frac{\alpha}{c_0^2}\Norm{\rho}^2
    &\leq
    \norm{\cI_1} + \norm{\cI_2}
    + \norm{\cI_3} + \norm{\cI_4}.
  \end{split}
\end{equation}
Solving the differential equation with an integrating
factor approach and integrating from $0$ to $T$ we see
\begin{equation}
  \frac{1}{2}
  \Norm{\rho(T)}^2 
  -
  \frac{1}{2}
  \Norm{\rho(0)}^2  
  \leq 
  \int_0^T
  \exp{\left(\frac{\alpha}{c_0^2}(s-T)\right)} 
  \left(
   \norm{\cI_1} + \norm{\cI_2}
    + \norm{\cI_3} + \norm{\cI_4}
    \right)
  \d s.
\end{equation}
Denote $t_* \in [0,T]$ to be the time such that
\begin{equation}
  \Norm{\rho(t_*)} = \max_{t\in[0,T]}\Norm{\rho(t)}
\end{equation}
we see that
\begin{equation}
  \frac{1}{2}\Norm{\rho(t_*)}^2
  -
  \frac{1}{2}\Norm{\rho(0)}^2
  \leq
  \int_0^{t_*}
  \exp{\left(\frac{\alpha}{c_0^2}(s-t_*)\right)} 
  \left(   \norm{\cI_1} + \norm{\cI_2}
    + \norm{\cI_3} + \norm{\cI_4}
    \right)
  \d s.
\end{equation}
It then follows that
\begin{equation}
  \begin{split}
    \frac{1}{2}\Norm{\rho(t_*)}^2 
    -
    \frac{1}{2}\Norm{\rho(0)}^2
    &\leq
    \int_0^{T}
    \exp{\left(\frac{\alpha}{c_0^2}(s-T)\right)} 
    \left(
    \norm{\cI_1} + \norm{\cI_2}
    + \norm{\cI_3} + \norm{\cI_4}
    \right)
    \d s.
    \\
    &\leq
    \sum_{n=1}^N \int_{t_{n-1}}^{t_n} 
    d(s)
    \left(
    \norm{\cI_1} + \norm{\cI_2}
    + \norm{\cI_3} + \norm{\cI_4}
    \right)
    \d s.
  \end{split}
\end{equation}
The terms $\cI_1$, $\cI_3$, $\cI_4$ are all
dealt with in a similar way, by using \CBS
and a maximum argument. For example
\begin{equation}
  \begin{split}
    \int_{t_{n-1}}^{t_n}
    d(s)
    \norm{\cI_1} \d s
    &\leq
    \int_{t_{n-1}}^{t_n}
    d(s)
    \norm{    \ltwop{\partial_t\epsi(s)}{\rho(s)}} \d s
    \\
    &\leq
    \int_{t_{n-1}}^{t_n}
    d(s)
    \Norm{\partial_t\epsi(s)}\Norm{\rho(s)} \d s
    \\
    &\leq
    \Norm{\rho(t_*)}
    \int_{t_{n-1}}^{t_n}
    d(s)
    \Norm{\partial_t\epsi(s)} \d s.
    \\
    &\leq
    \Norm{\rho(t_*)}
    \int_{t_{n-1}}^{t_n}
    d(s)
    \Norm{\epsilon^n - \epsilon^{n-1}} \d s.
    \\
    &\leq
    \Norm{\rho(t_*)}
    \int_{t_{n-1}}^{t_n} 
    d(s)
    \eta_n \d s.
  \end{split}
\end{equation}
The term $\cI_2$ that will eventually yield a time error
indicator requires a little more care.
\begin{equation}
  \begin{split}
    \int_{t_{n-1}}^{t_n}
    d(s)
    \norm{\cI_2} \d s
    &=
    \int_{t_{n-1}}^{t_n}
    d(s)
    \norm{\abil{\w(s) - \w^n}{\rho(s)}} \d s
    \\
    &=
    \int_{t_{n-1}}^{t_n}
    d(s)
    \norm{\abil{l_{n-1}(s)\rec {n-1}U^{n-1}
        + l_n(s)\rec n U^n - \rec n U^n}{\rho(s)}} \d s
    \\
    &=
    \int_{t_{n-1}}^{t_n}
    d(s)
    l_{n-1}(s)
    \norm{\abil{\rec {n-1}U^{n-1}
        - \rec n U^n}{\rho(s)}} \d s
    \\
    &=
    \int_{t_{n-1}}^{t_n}
    d(s)
    l_{n-1}(s)
    \norm{\ltwop{A^{n-1} U^{n-1} - A^n U^n}{\rho(s)}} \d s
    \\
    &\leq
    \Norm{\rho_*}
    \int_{t_{n-1}}^{t_n}
    d(s)
    \Norm{A^{n-1} U^{n-1} - A^n U^n} \d s
    \\
    &\leq
    \Norm{\rho_*}
    \int_{t_{n-1}}^{t_n}
    d(s)
    \indtime n \d s
  \end{split}
\end{equation}
Combining the results together we see
\begin{equation}
  \Norm{\rho(t_*)}^2 
  \leq
  \Norm{\rho(0)}^2 
  + 
  2 \Norm{\rho(t_*)}
  \sum_{n=1}^N
  \int_{t_{n-1}}^{t_n}
  d(s)
  \left(
  \eta_n + \theta_n + \beta_n + \gamma_n
  \right) \d s.
\end{equation}
Making use of the $\leb{2}(\W)$ simplification rule \citep[\S
  3.8]{lakkis-makridakis:06} it follows that
\begin{equation}
  \Norm{\rho(t_*)}
  \leq
  \Norm{\rho(0)} 
  + 
  2
  \sum_{n=1}^N
  \int_{t_{n-1}}^{t_n}
  d(s)
  \left(
  \eta_n + \theta_n + \beta_n + \gamma_n
  \right) \d s,
\end{equation}
which yields the desired result.
\end{proof}
\section{Numerical comparison of duality and energy, via spatial residuals}
\label{sec:sample-application-residuals}
We close the paper with a sample application of the ``abstract''
\aposteriori estimates derived in
\S\S\ref{sec:energy}--\ref{sec:duality.fully}.  

In particular, we summarise next numerical experiments to test the
asymptotic behaviour of the estimators given in Theorem
\ref{the:general.aposteriori.estimate}, Corollary
\ref{cor:fully-discrete-estimate} and Theorem
\ref{the:general.energy.aposteriori.estimate}. The C code used for
these computational experiments is based on the library \alberta
\citep{alberta}. To make the effects of numerical quadrature
negligible we choose the quadrature formula such that it is exact on
polynomials of degree 17 and less.
\subsection{Residuals}
Since our aim is to compare the numerical performance of the duality
and the energy based estimates, which differ mostly in their
time-accumulation and time-estimation aspects, we use the same type of
spatial indicators given by the residual estimators function
introduced in Lemma~\ref{lem:elliptic.aposteriori.estimates} and refer
to \citet{lakkis-makridakis:06} or \citet{LakkisPryer:10} for more
details.

The residuals constitute the building blocks of the \aposteriori
estimators used in our computer experiments.  We associate with
equations \eqref{eqn:continuous.heat} and \eqref{eqn:pointwise.form}
two residual functions: the {\em inner residual} is defined as
\begin{equation}
\label{eqn:inner.residual.function}
\begin{split}
  R^0 &:= \pwA  U^0-A^0U^0,\\
  R^n &:= \pwA \un-\an\un = \pwA\un-\projf n+\discn\un,\text{ for
  }\rangefromto n1N,
\end{split}
\end{equation}
and the {\em jump residual} which is defined as
\begin{equation}
\label{eqn:jump.residual.function}
J^n := J[\un]=\jump{\grad U^n}.
\end{equation}
With definition \secref{obs:representation} in mind, the inner
residual terms can be written explicitly as
\begin{equation}
  \ltwop{R^n}\phi
  =\sum_{K\in\cT_n}
  \ltwop{\contopera v-\projf n
    +\frac{U^n-\lproj n U^{n-1}}{\taun}}\phi _K.
\end{equation}

We can now introduce, for $\rangefromto n0N$, the {\em elliptic
reconstruction error indicators}
  \begin{gather}
    \label{eqn:lil2.indrec}
    \indrec n
    :=C_{6,2}\Norm{\hn^2 R^n }+C_{10,2}\Norm{\hn^{3/2} J^n }_{\Sigma_n}
    ,
    \\
    \intertext{and, for $\rangefromto n1N$, the {\em space error indicator}}
    \label{eqn:lil2.indspace}
    \indspace_n
    :=C_{6,2}\Norm{\hathn^2\bdisc R^n}
    +C_{10,2}\Norm{\hathn^{3/2}\bdisc  J^n}_
    {\hat\Sigma_n}
    +C_{14,2}\Norm{\hathn^{3/2}\bdisc  J^n}_
    {\check\Sigma_n\Tolto\hat\Sigma_n}
    .
\end{gather}
\subsection{The benchmark problem}
We take $\geovec A = -\geovec I$ such that the parabolic
problem~(\ref{eqn:continuous.heat}) coincides with the heat
equation. We tune data functions $f$ and $u_0$ of this parabolic problem
so that its exact solution $u$ is given by
\begin{equation}
  \label{eqn:Problem}
  u(\vec x, t) 
  =
  \sin{\left(\kappa \pi t\right)}\exp{\left(-10\norm{\vec x}^2\right)},
\end{equation}
with $\kappa\in\naturals$. We fix $d=2$ and take $\W = [-1,
  1]\times[-1, 1]$.
\begin{Def}[experimental order of convergence]
  Given two sequences $a(i)$ and $h(i)$, $i=0,\ldots,N$ we define
  the experimental order of convergence (\EOC) to be:
  \begin{equation}
    \label{eqn:def:EOC}
    \EOC(a,h;i) = \frac{ \ln(a(i+1)/a(i)) }{ \ln(h(i+1)/h(i)) }.
  \end{equation}
\end{Def}
\begin{Def}[effectivity index and its inverse]
  The main tool deciding the quality of an estimator is the
  effectivity index (\EI) which is the ratio of the error and the
  estimator, \ie using the estimators from the duality-based
  Theorem \ref{the:general.aposteriori.estimate} at time $t_m$, for
  some $\rangefromto m1N$
  \begin{equation}
    \EI(t_m) = 
    \frac{
      \psqrt{\smash{
        \sum_{n=1}^m b_n \theta_n^2
      }}
      +
      \psqrt{\smash{
        \sum_{n=0}^{m-1} a_n\epsilon_n^2
      }}
      +
      \eta_m
    }{
      \Norm{U-u}_{\leb{\infty}(0,t_m;\leb{2}(\W))}
    },
  \end{equation}
  using the results of Corollary \ref{cor:fully-discrete-estimate}
  \begin{equation}
      \EI(t_m) = 
      \frac{
        \max_{n\in[0:m]} \theta_n
	+
        \max_{n\in[0:m-1]}\epsilon_n
        +
        \eta_m
      }
      {\Norm{U-u}_{\leb{\infty}(0,t_m;\leb{2}(\W))}}
    \end{equation}
    and for the estimator associated with the energy-estimator Theorem
    \ref{the:general.energy.aposteriori.estimate}
    \begin{equation}
      \EI(t_m) = 
      \frac{
        \sum_{n=1}^m \int_{t_{n-1}}^{t_n} d(s) (\theta_n + \eta_n)
        +
        \max_{n\in[0:m]}\epsilon_n
      }
           {\Norm{U-u}_{\leb{\infty}(0,t_m;\leb{2}(\W))}}.
    \end{equation}
    Since it is much
    easier to visualise we will be computing the \emph{inverse
      effectivity index}, $1/\EI(t_m)$.
  \end{Def}
\subsection{Comparing duality estimates with energy estimates}
The second main objective of this research was to compare the
numerical results associated with the duality-based estimator
(Th. \ref{the:general.aposteriori.estimate}) to those of the
energy-based estimator
(Them. \ref{the:general.energy.aposteriori.estimate}).  As already
observed in Figure~\ref{fig:time-accumulation-coefficients}, we expect
the energy-based estimators to perform better due to a
time-accumulation of the estimator which, due to the exponential tail
of the integration weight, is more consistent with the
$\leb\infty(0,T;\leb2(\W))$ norm.  Extensive numerical experimentation
lead to the similar conclusions: energy-based estimators perform
slightly better than duality-based ones for this norm for short-time
integration and much better for long-time integration.  For space
reasons we give here only the results on one example, illustrating
this point.

The benchmark problem (\ref{eqn:Problem}) has been chosen such that
our results can be compared with those in
\citet{lakkis-makridakis:06}. The initial condition is zero, the
boundary values are not exactly zero but negligible, hence little
interpolation error is committed here; however some care must be taken
dealing with these small numbers.  The model problem
\eqref{eqn:Problem} is approximated on a stationary mesh in time,
hence $\fes{j-1}=\fes{j}$ for all $j=1,\dotsc,N$, for two values of
the \emph{time-oscillation factor} $\kappa=1$ and $8$.  Our results
are valid for any polynomial order for the spatial finite elements,
but we report results only for $\poly1$ elements. In order to
emphasize the time-estimator, we take $\tau\equiv h$ in all these
experiments.

All estimators appearing in
Theorems~\ref{the:general.aposteriori.estimate},~\ref{cor:fully-discrete-estimate}
and \ref{the:general.energy.aposteriori.estimate} are computed except
for $\indcoarse n$ and $\inddata n$.  The mesh-change estimator
$\indcoarse n = 0$ in our tests here, because nowhere is the
triangulation $\cT^{n}$ a coarsening of $\cT^{n-1}$; for examples with
$\gamma_n\neq0$ we refer to \citet{LakkisPryer:10}.  We do not track
the data error indicator $\beta_n$ either, since it can be shown to be
of higher order in our case, due to regularity of $f$.

In Figures \ref{Fig:duality-P1h}--\ref{Fig:duality-P1h2} we plot
convergence of error and estimators derived via duality in Theorem
\ref{the:general.aposteriori.estimate} and Corollary
\ref{cor:fully-discrete-estimate}. In Figures \ref{Fig:energy-P1} we
are report the corresponding results for the estimator derived via the
energy argument in
Theorem~\ref{the:general.energy.aposteriori.estimate}.  The results,
commented in Figure~\ref{Fig:duality-P1h}--\ref{Fig:energy-P1}'s
captions, confirm the energy-based estimator's superiority, slightly for short
times ($k=1$), and clearly for long integration times ($k=8$).
\providecommand{\figwidth}{0.8\textwidth} \newcommand{\figscale}{1}
  \begin{figure}[ht]
    \caption{\label{Fig:duality-P1h}
      Convergence of error and the duality estimators of
      \S\ref{sec:duality.fully} for \eqref{eqn:Problem} with low
      ``time-oscillation factor'' $\kappa = 1$ and $\poly1$ elements
      on uniform meshes and timestep with $\tau=0.05 h$ and $h(i) =
      2^{-i}$, $\rangefromto i5{11}$.  We plot all quantities as
      functions of (PDE) time. Rates for each error/indicator can be
      read from the experimental order of convergence (EOC) of the
      associated part of the estimator together with its value on a
      logarithmic scale.  The colour/grey scale is such that dark is
      finest and light is coarsest.  Since the benchmark problem has
      no initial error, we employ the inverse effectivity index.
    }
    \begin{center}
      \subfigure[][{
          \label{subfig:duality-P1h-L2-accumulation-prob1}
          We test duality estimators from Theorem \ref{the:general.aposteriori.estimate}.}]{
        \includegraphics[scale=\figscale,width=\figwidth]
                        {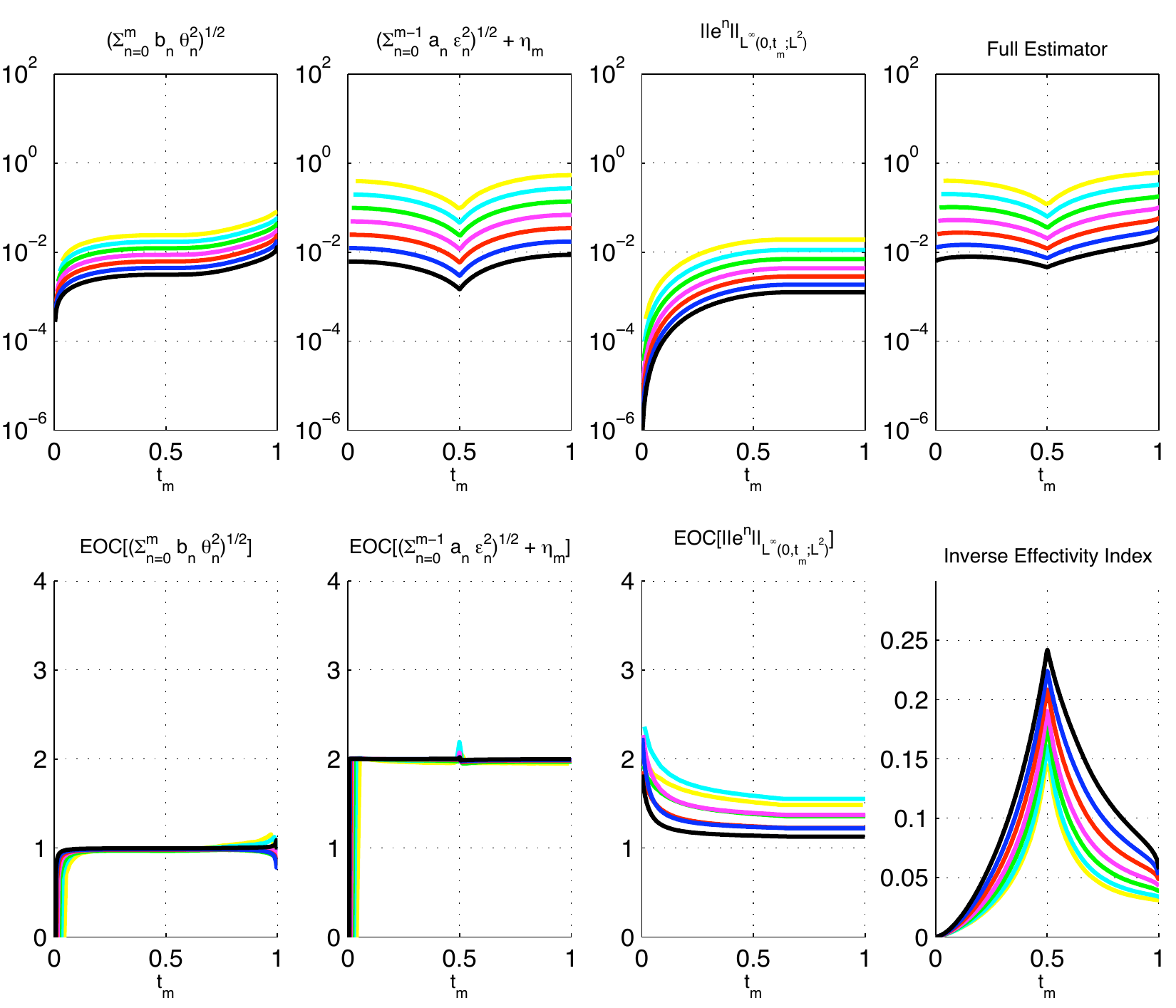}
      }
      \subfigure[][{Here we study the estimator from Corollary
          \ref{cor:fully-discrete-estimate}.}]{
        \label{subfig:duality-P1h-max-accumulation-prob1}
      \includegraphics[scale=\figscale,width=\figwidth]
                      {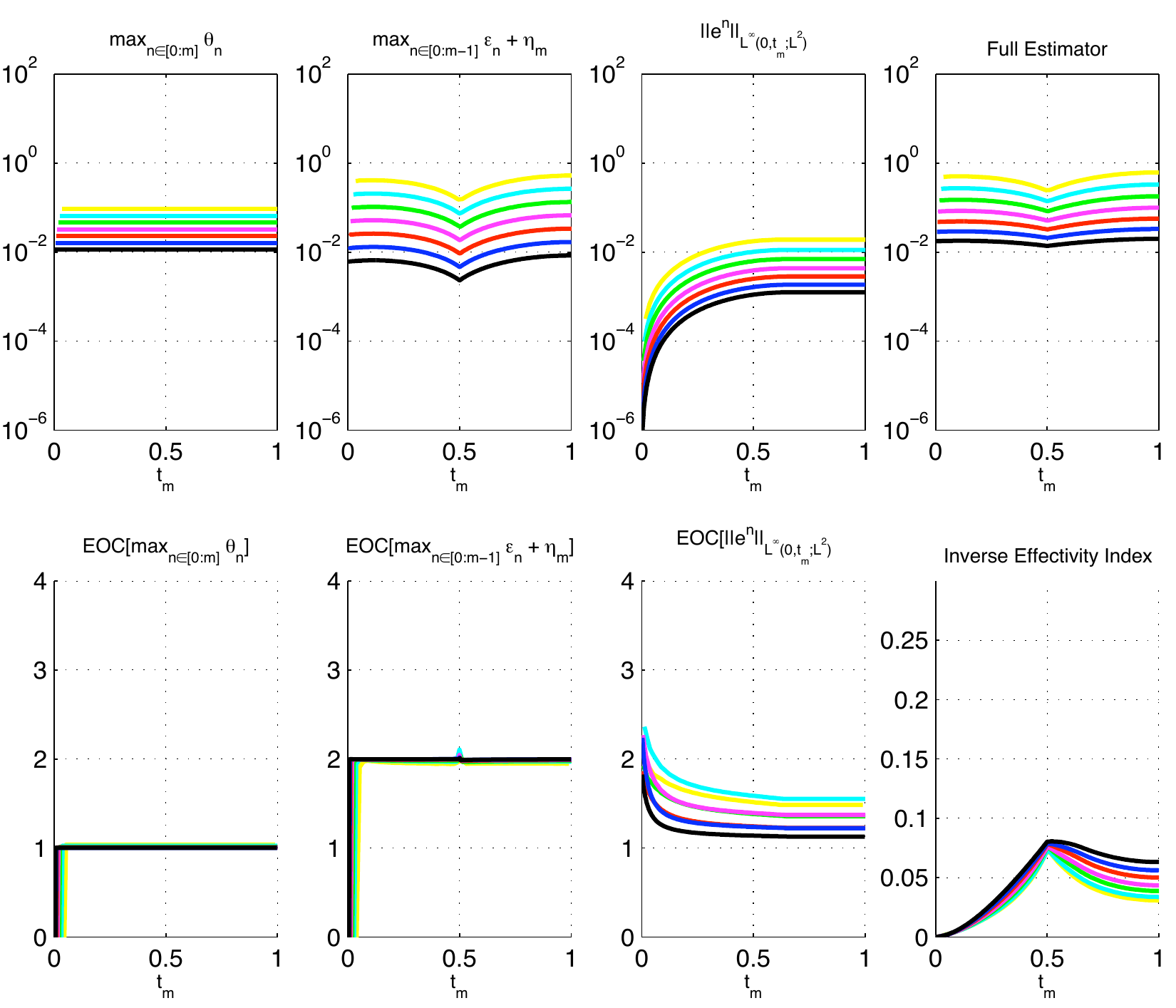}
                        
      }
    \end{center}
  \end{figure}

  \begin{figure}[ht]
    \caption{\label{Fig:duality-P1h2}
      Convergence of error and the duality estimators for
      \eqref{eqn:Problem} with time-oscillation factor $\kappa=8$ and
      the other parameters as in Figure~\ref{Fig:duality-P1h}.  We
      plot all quantities as functions of (PDE) time. Rates for each
      error/indicator can be read from the experimental order of
      convergence (EOC) of the associated part of the estimator
      together with its value on a logarithmic scale.  The colour/grey
      scale is such that dark is finest and light is coarsest.  Since
      the benchmark problem has no initial error, we employ the
      inverse effectivity index.  Note that due to the faster
      oscillation, how a discrepancy between the error's
      time-accumulation (or rather the lack thereof) and the
      estimator's builds up quickly and produces oscillations in the effectivity index.
    }
    \begin{center}
      \subfigure[][{Here we study the estimator from
          Theorem \ref{the:general.aposteriori.estimate}.}]{
        \label{subfig:duality-P1h2-L2-accumulation-prob2}
      \includegraphics[scale=\figscale,width=\figwidth]
                      {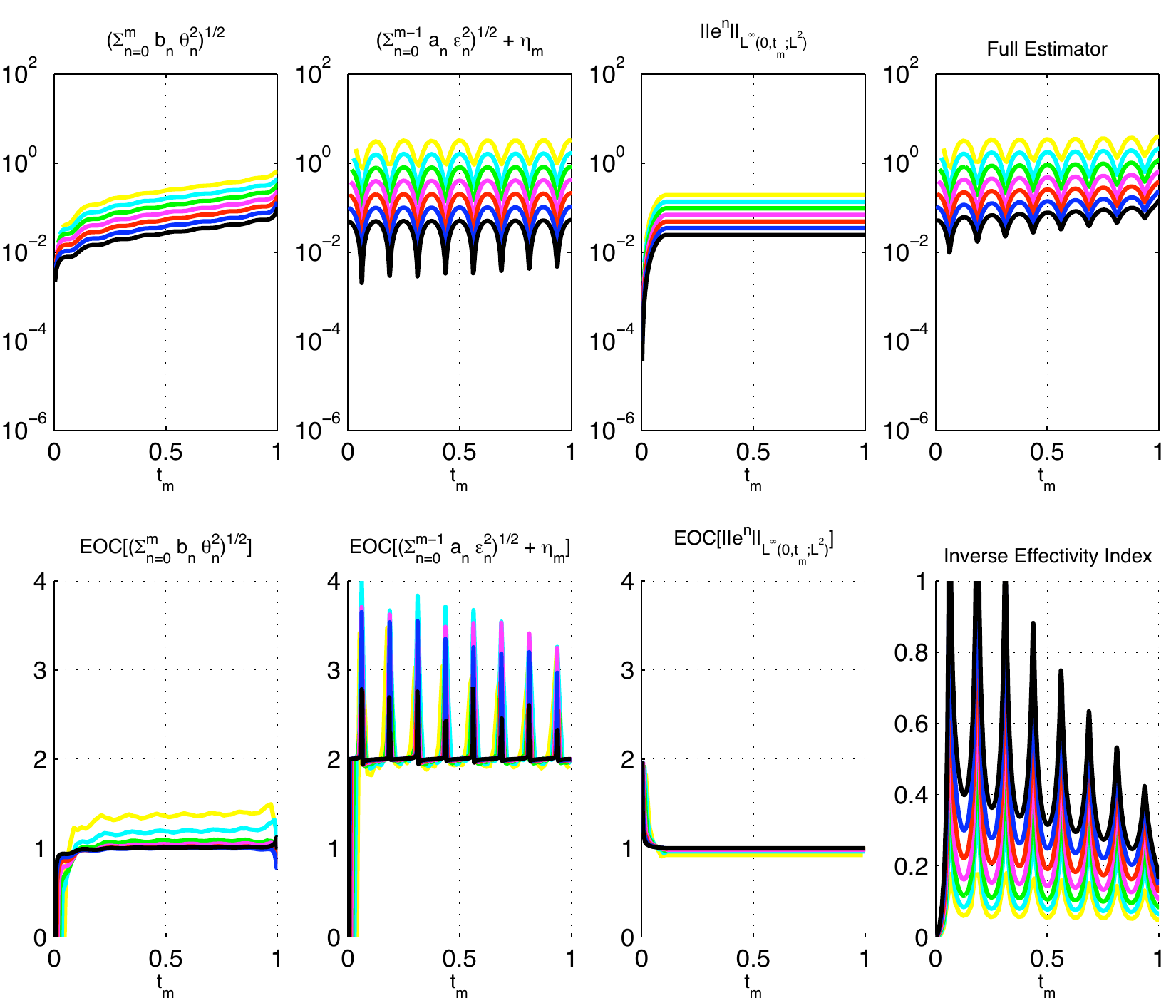} 
      }
      \subfigure[][{Here we study the estimator from Corollary
          \ref{cor:fully-discrete-estimate}.}]{
        \label{subfig:duality-P1h2-max-accumulation-prob2}
      \includegraphics[scale=\figscale,width=\figwidth]
                      {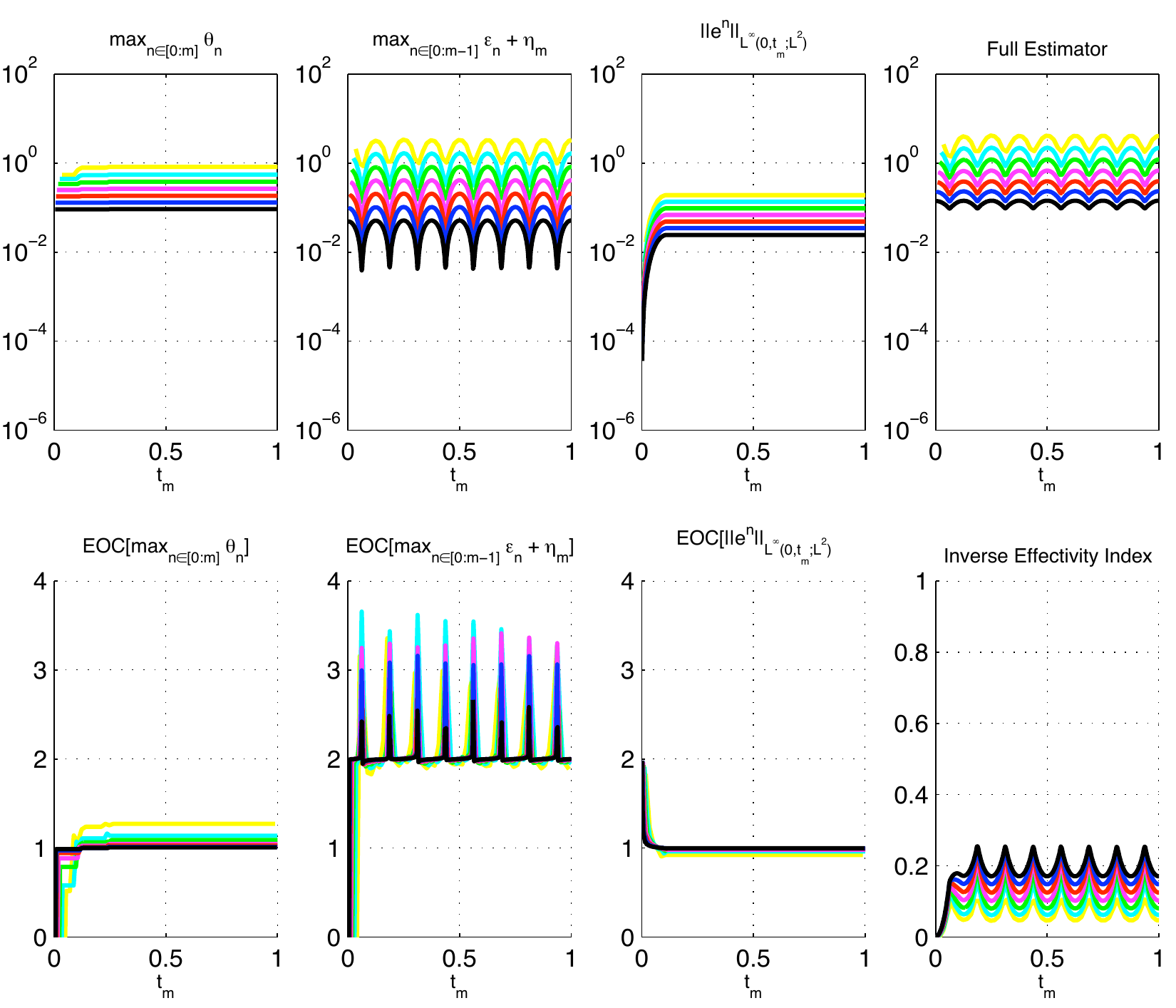}
      }
    \end{center}
  \end{figure}
  \begin{figure}[ht]
    \caption{
      \label{Fig:energy-P1} 
      Convergence for of the energy-based estimator from Theorem
      \ref{the:general.energy.aposteriori.estimate}, taking a uniform
      timestep $\tau = 0.05 h$. We couple $\tau$ to appropriate powers
      of $h$ and with $h(i) = 2^{-i}$, $\rangefromto i5{11}$. In both
      cases we are using $\poly1$ elements. We plot the \EOC of the
      associated part of the estimator together with its value on a
      logarithmic scale.
    }
    \begin{center}
      \subfigure[][{
          Problem \eqref{eqn:Problem} with time-oscillation
          factor $\kappa = 1$, to be compared with
          Figure~\ref{Fig:duality-P1h}.  Convergence rates are clearly
          similar, but the effectivity index is much more stable (and
          smooth) with respect to time in this case, as expected given
          the better accumulation of the estimator in time (cf.
          Figure~\ref{fig:time-accumulation-coefficients}).
      }]{
        \label{subfig:energy-P1h-prob1}        
        \includegraphics[scale=\figscale,width=\figwidth]
        {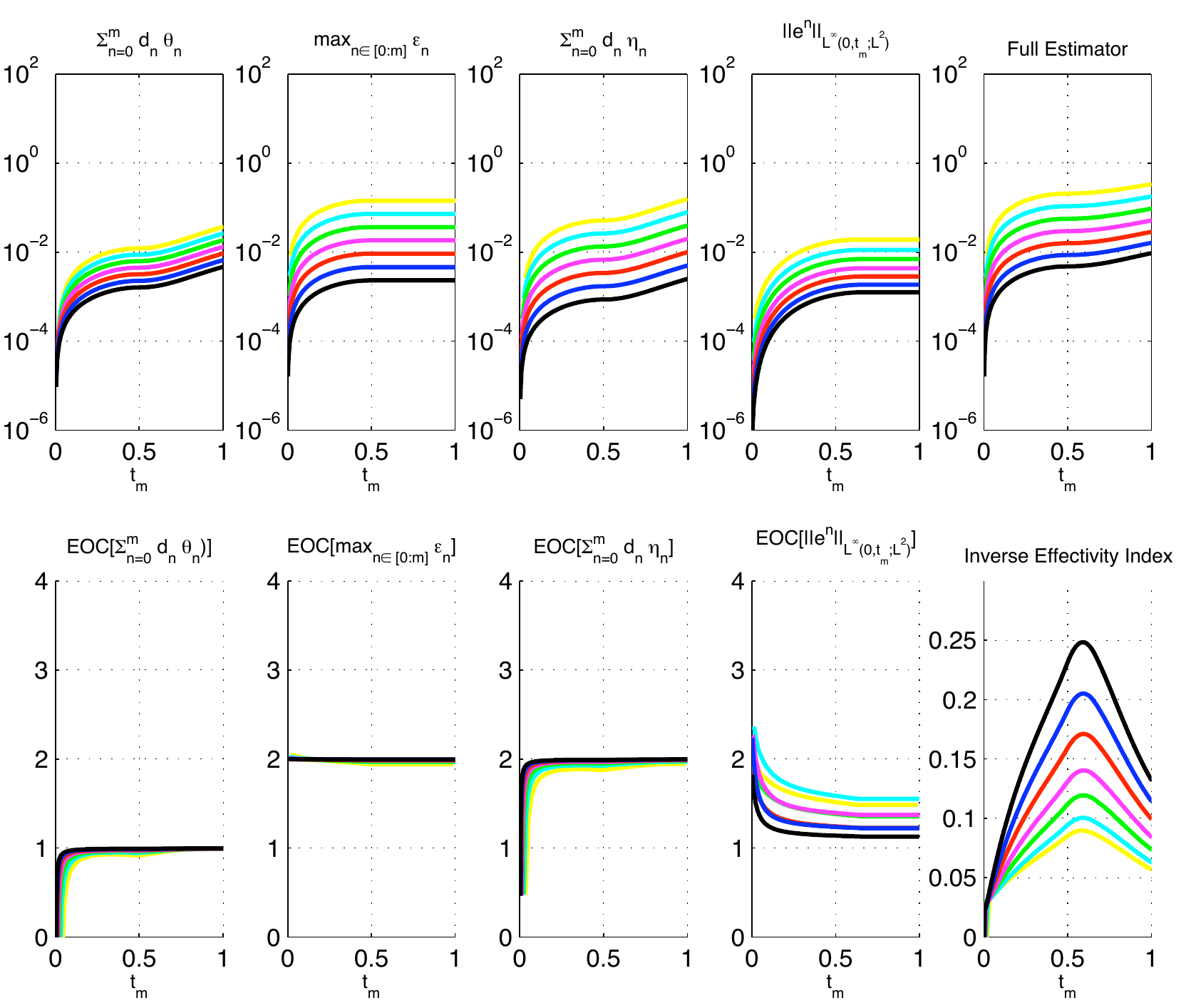}
      } 
      \subfigure[][{
          Problem \eqref{eqn:Problem} with time-oscillation factor
          $\kappa = 8$, to be compared with
          Figure~\ref{Fig:duality-P1h2}.  Here the gain in effectivity
          index, from using the energy instead of duality estimators is even more
          dramatic.
          }]{
        \label{subfig:energy-P1h2-prob2}        
        \includegraphics[scale=\figscale,width=\figwidth]
                        {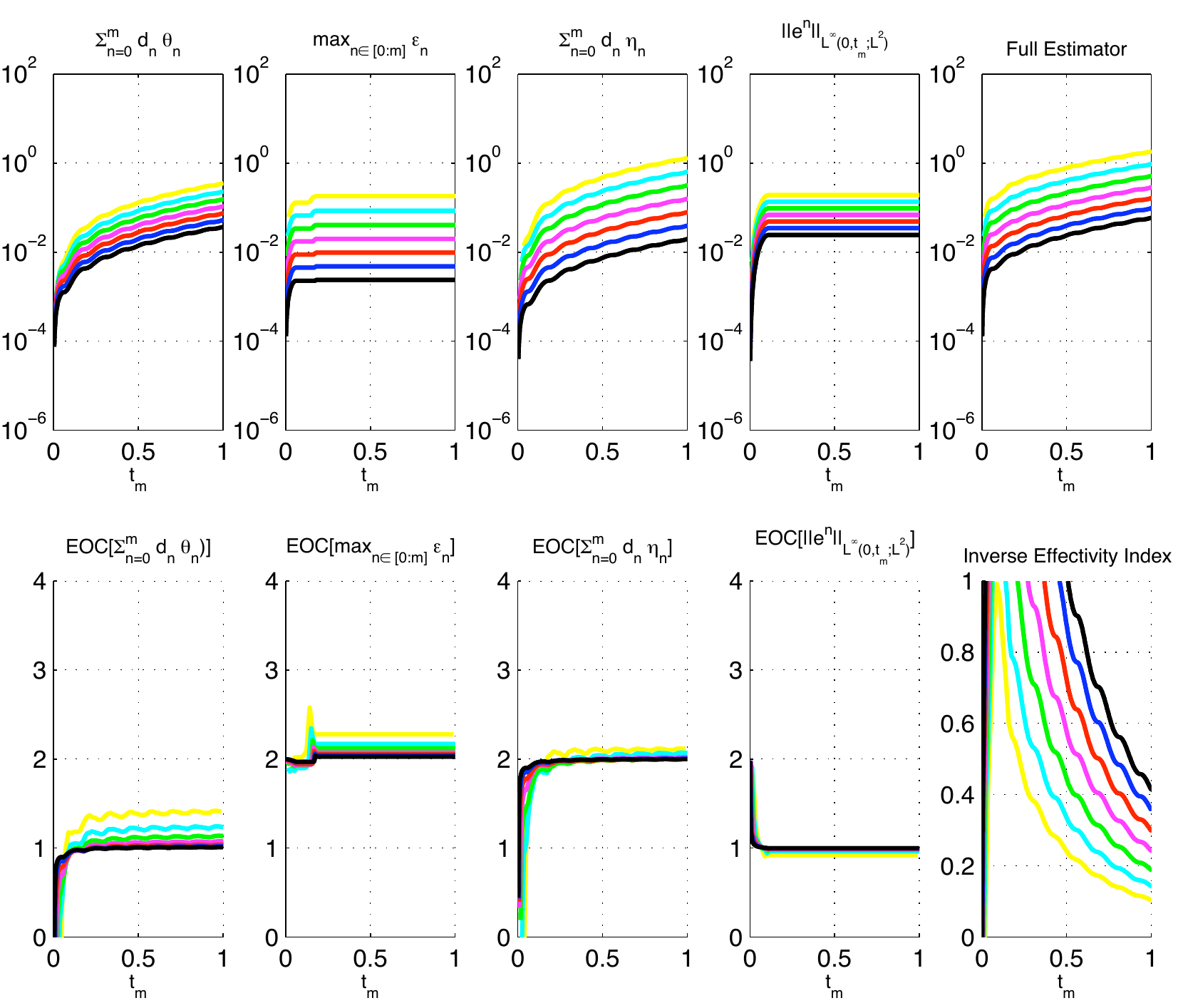}
      }
    \end{center}
  \end{figure}
\end{document}